\documentclass[arxiv,reqno,twoside,a4paper,12pt]{amsart}

\usepackage{tikz}
\usepackage{latexsym}
\usepackage{amssymb}
\usepackage{amsmath}
\usepackage{amsfonts}
\usepackage{amsthm}
\usepackage{thmtools}

\usepackage{nameref}
\usepackage[colorlinks=true,linkcolor=blue,citecolor=blue]{hyperref}
\usepackage[driver=pdftex,heightrounded=true,centering]{geometry}
\usepackage{cleveref}
\usepackage{chngcntr}
\usepackage{cases}

\usepackage[scaled]{helvet} 
\usepackage{courier} 
\usepackage[mathbf]{euler}

\counterwithin{equation}{section}

\declaretheorem[name=Theorem,numberwithin=section]{theo}

\declaretheorem[name=Lemma,sibling=theo]{lem}
\declaretheorem[name=Corollary,sibling=theo]{cor}
\declaretheorem[name=Remark,sibling=theo]{rem}

\crefname{theo}{Theorem}{Theorems}
\crefname{pro}{Proposition}{Propositions}
\crefname{lem}{Lemma}{Lemmas}
\crefname{cor}{Corollary}{Corollaries}
\crefname{rem}{Remark}{Remarks}
\crefname{rems}{Remarks}{Remarks}
\crefname{defi}{Definition}{Definitions}
\crefname{example}{Example}{Examples}

\newcommand{\wt}{\widetilde}

\newcommand{\ov}{\overline}

\renewcommand{\d}{\delta}
\renewcommand{\l}{\lambda}

\newcommand{\ve}{\varepsilon}

\newcommand{\g}{\gamma}
\newcommand{\G}{\Gamma}

\renewcommand{\r}{\mathbb{R}}
\newcommand{\N}{\mathbb{N}}

\renewcommand{\sec}{\operatorname{sec}}

\DeclareMathOperator{\Rm}{Rm}
\DeclareMathOperator{\Ric}{Ric}
\renewcommand{\le}{\leq}
\newcommand{\p}{\partial}
\renewcommand{\a}{\alpha}
\renewcommand{\b}{\beta}

\newcommand{\n}{\nabla}

\title{Ricci de Turck flow on incomplete manifolds}
\author{Tobias Marxen}
\author{Boris Vertman}
\date{\today}

\numberwithin{equation}{section}

\begin{document}

\begin{abstract}
In this paper we construct a Ricci de Turck flow on any incomplete Riemannian manifold with bounded curvature.
The central property of the flow is that it stays uniformly equivalent to the initial incomplete Riemannian metric, and in that sense 
preserves any given initial singularity structure. Together with the corresponding result by Shi for complete manifolds \cite{shi}, 
this gives that any (complete or incomplete) manifold of bounded curvature can be evolved by the Ricci de Turck flow for a short time. 
\end{abstract}

\maketitle
\tableofcontents

\section{Introduction and statement of the main result}\label{intro}

\bigskip

Consider an $n$-dimensional, smooth and possibly incomplete Riemannian manifold $(M,\wt{g})$.
We denote the corresponding Riemannian curvature tensor by $\wt{\Rm}$ and its pointwise norm with respect to 
$\wt{g}$ by $|\wt{\Rm}|$. The Ricci de Turck flow of $(M,\wt{g})$ is a smooth family $g(t), t\in [0,T],$
of Riemannian metrics on $M$, solving the initial value problem
\begin{equation}\label{Ricci de Turck system on M incomplete}
\frac{\p}{\p t} g_{ij}(t) = -2 \Ric_{ij}(t) + \n_i V_j(t) + \n_j V_i(t), \quad
g(0) = \wt{g}.
\end{equation}
where $V^i(t) = g(t)^{jk}(\G^i_{jk}(g(t)) - \G^i_{jk}(\wt{g}))$ is the de Turck vector field
defined\footnote{We employ the Einstein summation convention.} in terms of Christoffel symbols $\G^i_{jk}$ for $g(t)$ and $\wt{g}$;
$(\Ric_{ij}(t))$ is the Ricci curvature tensor and $\n$ the covariant derivative of $g(t)$.
Our main theorem is then as follows.

\begin{theo}\label{main theo}
Assume $|\wt{\Rm}|^2 \le k_0$ for some positive constant $k_0>0$. 
Then there exists $T(n,k_0) > 0$, depending only on $n$ and $k_0$, such that the initial value problem
\eqref{Ricci de Turck system on M incomplete} has a smooth solution $g(t)$ for $t\in [0,T(n,k_0)]$.
Furthermore, for any $\delta > 0$ there exists $0 < T(n,k_0,\delta) \le T(n,k_0)$ depending only on $n, k_0$ and $\delta$, such that
\begin{equation} \label{unif equiv delta on M incomplete}
(1-\delta) \wt{g}(x) \le g(x,t) \le (1+\delta) \wt{g}(x),
\end{equation}
for all $(x,t) \in M \times [0,T(n,k_0,\delta)]$. Moreover, if we assume that for all $m \ge 1$ there exists a constant $C_m > 0$, such that for all $x \in M$, $0 < \rho \le 1$
\[ |\wt{\n}^m \wt{\Rm}|(x) \le \frac{C}{\rho^m}  \] 
whenever $B(x,\rho-r)$ is relatively compact for all $r > 0$, then there exist constants $C' > 0$, $C'_m > 0$, such that for all $x \in M$, $t \in [0,T]$, $0 < \rho \le 1$
\[ |\wt{\n}^m g|(x,t) \le \frac{C'_m}{\rho^m}, \quad |\Rm|(x,t) \le \frac{C'}{\rho^2}, \quad |\n^m \Rm|(x,t) \le \frac{C'}{\rho^{m+2}} \]
whenever $B(x,\rho-r)$ is relatively compact for all $r > 0$.	
\end{theo}
	
\begin{rem}
The condition that $B(x,\rho-r)$ is relatively compact in $M$ for all $r > 0$ is an intrinsic way to express the distance of a point 
$x \in M$ to the singular strata of $M$. It means that this distance is larger or equal to $\rho$. 	
\end{rem}

We should point out that short-time existence and further properties of a Ricci de Turck flow on incomplete manifolds has already
been established in the special case of manifolds with conical or more generally wedge singularities in varying dimensions
in \cite{MRS}, \cite{BV}, \cite{v}, \cite{Klaus-Vertman} and \cite{Yin:RFO}, to name a few. These references deal with the
flow that stays uniformly equivalent to the initial metric and hence preserves the initial singularity. Due to non-uniqueness of
the flow in the singular setting, there exist solutions that are instantaneously complete, cf. \cite{Topping2}, 
as well as solutions that smooth out the singularity, cf. \cite{MS}. \medskip

The main novelty of the present paper is the assertion that such a Ricci de Turck flow, preserving the initial 
singularity structure, exists on any \emph{arbitrary} incomplete manifold of bounded curvature. This includes, 
but is not restricted to, for instance incomplete $3$-dimensional manifolds with isolated conical singularities, 
where the singularity is a Ricci-flat cone in first approximation.
In this setting we also establish explicit estimates for arbitrary higher derivatives of the metric 
and of the Riemann curvature tensor along the flow. We conjecture that this flow coincides with the flows 
studied in our previous works such as e.g. in \cite{v}. \medskip

Our paper is structured as follows. In \S \ref{Review of Shis local existence theorem} we review the argument of 
Shi \cite{shi}, which proves short time existence of Ricci de Turck flow for complete manifolds of bounded curvature. We break down the argument
to those points where completeness of the manifold is used. In the subsequent \S \ref{curvature-section},
\S \ref{curvature-section-second-order} and \S \ref{curvature-section-higher-order} we establish a priori estimates 
for the first, second and higher derivatives of the metric along the flow. \S \ref{curvature-section-second-order} and \ref{curvature-section-higher-order} also contain a priori estimates for the Riemann curvature tensor. In the final \S \ref{main-section} we adapt the
argument of \S \ref{Review of Shis local existence theorem} in order to establish the corresponding result for incomplete manifolds of bounded curvature as well. 
\medskip

\textit{Notation:} Let us fix the notation for the discussion below. Let $g(t)$, $t \in [0,T]$ be a family of Riemannian metrics on an 
incomplete manifold $M$. We denote by $\n$ and $\G$ the covariant derivative and the Christoffel symbols with respect to $g(t)$. 
$\Rm$, $\Ric$ and $R$ denote the Riemann curvature tensor, the Ricci tensor and the scalar curvature of $g(t)$, 
respectively.  \medskip

Let $\wt{g}$ be the initial Riemannian metric on $M$. Quantities with respect to $\wt{g}$ are marked with an
upper tilde. For example we write $\wt{\n}$ for the covariant derivative with respect to $\wt{g}$. There are the following exceptions to this rule:
We denote by $B(x,r)$ the open ball with radius $r > 0$ and centre $x \in M$, and we write $B(A,r) := \{x \in M : d_{\wt{g}}(x,A) < r\}$ 
for the $r$-neighborhood of a given subset $A \subset M$, both with respect to the metric $\wt{g}$. The norm $|\cdot|$ will always be 
with respect to $\wt{g}$. We write $d_{\wt{g}}$ for the distance function induced by $\wt{g}$.

\section{Review of Shi's local existence theorem} \label{Review of Shis local existence theorem}

In this section we review results and proofs from Shi \cite{shi} in the complete setting.
Shi established the following short-time existence result for the Ricci de Turck flow starting at 
complete manifolds with bounded curvature. Within this section, $(M,\wt{g})$ is always understood
to be a complete $n$-dimensional Riemannian manifold of bounded curvature. 

\begin{theo}[\cite{shi}, Theorems 4.3, 2.5]\label{Shi Theorems}
Assume $|\wt{\Rm}|^2 \le k_0$ for some positive constant $k_0 > 0$. 
Then there exists $T(n,k_0) > 0$ depending only on $n$ and $k_0$, such that the initial value problem
\eqref{Ricci de Turck system on M incomplete} has a smooth solution $g(t)$.
Moreover, for any $\delta > 0$ there exists $0 < T(n,k_0,\delta) \le T(n,k_0)$ depending only on $n, k_0$ and $\delta$, such that
\begin{equation} \label{unif equiv delta on M complete}
(1-\delta) \wt{g}(x) \le g(x,t) \le (1+\delta) \wt{g}(x)
\end{equation}
for all $(x,t) \in M \times [0,T(n,k_0,\delta)]$.
\end{theo}

\begin{rem}\label{lower-injectivity}
We emphasize that the lower bound on the injectivity radius 
does not enter in the definition of the time bounds $T(n,k_0), T(n,k_0,\delta) > 0$.
Indeed, the local existence result still holds on complete manifolds without a positive
lower bound on the injectivity radius. An obvious instance are manifolds with hyperbolic 
cusps, where Theorem \ref{Shi Theorems} still holds despite the injectivity radius tending to zero at the cusp.
\end{rem}

The proof of this theorem is based on three main steps. The first is an a priori estimate for the 
Ricci de Turck flow on a relatively compact domain $D \subset M$ with Dirichlet boundary conditions.
\begin{theo}[\cite{shi}, Theorem 2.5]\label{Shi apriori}
Let $D \subset M$ be a relatively compact domain, whose boundary $\p D$ is an $(n-1)$-dimensional, 
smooth, compact submanifold. Let $g(x,t)$, $t \in [0,T]$ be a solution of the initial boundary value problem
\begin{equation} \label{Dirichlet Ricci}
\begin{split}
\frac{\p}{\p t} g_{ij}(x,t) = (-2 \Ric_{ij} + \n_i V_j + \n_j V_i)(x,t)&, \quad  (x,t) \in D \times [0,T],  \\
g(x,t) = \wt{g}(x)&, \quad (x,t) \in \partial D \times [0,T], \\
g(x,0) = \wt{g}(x)&, \quad x \in D.
\end{split}
\end{equation}
where $V^i = g^{jk}(\G^i_{jk} - \wt{\G}^i_{jk})$ is the de Turck vector field. Then for any $\delta > 0$ 
there exists $T(n,k_0,\delta) > 0$ depending only on $n, k_0$ and $\delta$, such that
\begin{equation} \label{unif equiv delta}
(1-\delta) \wt{g}(x) \le g(x,t) \le (1+\delta) \wt{g}(x)
\end{equation}
for all $(x,t) \in M \times [0,\min \, \{T(n,k_0,\delta),T\}]$.
\end{theo}

\begin{proof}[Proof outline]
Shi controls the eigenvalues $\l_k(x,t)$ of $g(x,t)$ with respect to $\wt{g}(x)$ (i.e. the eigenvalues of 
$g(x,t)$ considered as a $(1,1)$-tensor using the metric $\wt{g}(x)$).  Shi defines a function 
\[ \varphi(x,t) = \sum_{k=1}^n \l_k(x,t)^{-m}, \]
where $m>0$ is sufficiently large only depending on $n$ and $\d$. Shi then shows that $\varphi$ satisfies a differential inequality
\[ 
\frac{\p \varphi}{\p t} \le g^{\a\b} \wt{\n}_\a \wt{\n}_\b \varphi + 2mn\sqrt{k_0} \cdot \varphi^{1+1/m}, 
\]
and applies the maximum principle to conclude $\varphi(x,t) \le 2n$ for all $(x,t) \in D \times [0,T]$.
This leads to the lower bound in \eqref{unif equiv delta}. The upper bound in \eqref{unif equiv delta} 
is then obtained by a similar procedure applied to the function 
\[ 
F(x,t) = \left( 1-\frac{1}{2n}\sum_{k=1}^n \l_k(x,t)^{\wt{m}}\right)^{-1},
\]
where $\wt{m} > 0$ is large enough and only depends on $n$ and $\d$.  
\end{proof}

The second step is the short-time existence of system \eqref{Dirichlet Ricci}.
\begin{theo}[\cite{shi}, Theorem 3.2] \label{Shi Theorem 3.2}
Let $D \subset M$ be a relatively compact domain, whose boundary $\p D$ is an $(n-1)$-dimensional, 
smooth, compact submanifold. Then there exists $T(n,k_0) > 0$ only depending on $n$ and $k_0$, such that 
the initial boundary value problem \eqref{Dirichlet Ricci} admits a unique smooth solution 
$g(x,t)$, $(x,t) \in D \times [0,T(n,k_0)]$.
\end{theo}

The third step are interior estimates for the derivatives of the metric, only depending on $\wt{g}$ and not on any specified boundary conditions.
\begin{lem}[\cite{shi}, Lemma 4.1]\label{Shi Lemma 4.1}
Fix $0 < \g, \d, T < \infty$, and let $g(x,t)$ be a smooth solution of 
the initial value problem
\begin{equation*}
\begin{split}
\frac{\p}{\p t} g_{ij}(x,t) = (-2 \Ric_{ij} + \n_i V_j + \n_j V_i)(x,t)&, \quad  (x,t) \in \, B(x_0,\g+\d) \times [0,T],  \\
g(x,0) = \wt{g}(x)&, \quad x \in \, B(x_0,\g+\d),
\end{split}
\end{equation*}
where $V^i = g^{jk}(\G^i_{jk} - \wt{\G}^i_{jk})$ is the de Turck vector field. Furthermore, assume that
\begin{equation*}
(1-\ve(n)) \wt{g}(x) \le g(x,t) \le (1+\ve(n)) \wt{g}(x)
\end{equation*}
for $\ve(n) > 0$ sufficiently small, only depending on $n$, and for all $(x,t) \in B(x_0,\g+\d) \times [0,T]$. Then there exists a positive constant $c(n,\g,\d,T,\wt{g}) > 0$, 
depending only on $n, \g,\d,T$ and $\wt{g}$, such that
\begin{equation*}
|\wt{\n} g(x,t)|^2 \le c(n,\g,\d,T,\wt{g})
\end{equation*}
for all $(x,t) \in B(x_0,\g+\frac{\d}{2}) \times [0,T]$.
\end{lem}

\begin{proof}[Proof outline] 
Shi defines for any $(x,t) \in B(x_0,\g+\d) \times [0,T]$ the function
\begin{equation} \label{phi}
\varphi(x,t) = a + \sum_{k=1}^n \l_k(x,t)^{m_0}, 
\end{equation} 
where $a, m_0$ are carefully chosen positive constants only depending on $n$, and $\l_k(x,t)$ are 
the eigenvalues of $g(x,t)$ with respect to $\wt{g}(x)$. Shi then shows that the function
\begin{equation} \label{psi}
\psi(x,t) := |\wt{\n}g|^2 \varphi(x,t) 
\end{equation}
satisfies
\begin{equation} \label{psi equation}
\frac{\p \psi}{\p t} \le g^{\a\b} \wt{\n}_\a \wt{\n}_\b \psi - \frac{1}{16} \psi^2 + c_0, 
\end{equation}
where $c_0 > 0$ is a constant only depending on $n$ and $\wt{g}$. Then Shi takes a nonincreasing cutoff function 
$\eta \in C^\infty(\r)$ such that $\eta \equiv 1$ on $(-\infty, 0]$, vanishing identically on $[1,\infty)$ as illustrated in 
Figure \ref{fig:CutOff1}. 

\begin{figure}[h]
	\begin{center}
		
		\begin{tikzpicture}[scale=1.3]
		\draw[->] (-0.2,0) -- (6,0);
		\draw[->] (2,-0.2) -- (2,2.2);
		
		\draw (1.8,2) node[anchor=north] {$1$};
		
		\draw (0,2) -- (2,2);
		\draw (2,2) .. controls (3.4,2) and (2.6,0) .. (4,0);
		\draw[dashed] (4,-0.2) -- (4,2.2);
		
		\draw (3.5,1) node {$\eta$};
		
		\draw (2,-0.5) node {$0$};
		\draw (4,-0.5) node {$1$};

		\end{tikzpicture}
		
		\caption{The cutoff function $\eta$.}
		\label{fig:CutOff1}
	\end{center}
\end{figure}
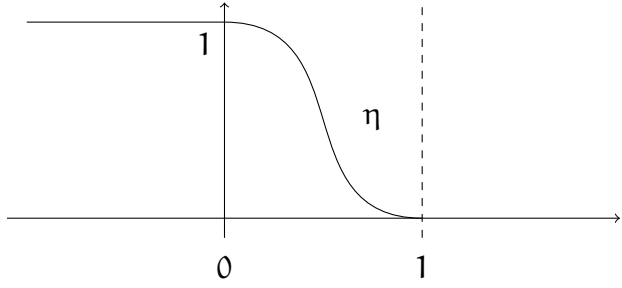

\noindent The crucial property of the function $\eta$ is the control on its derivatives
\begin{equation} \label{eta properties}
  |\eta''(x)| \le 8, \quad
  |\eta'(x)|^2 \le 16 \eta(x), \quad \textup{for any $x \in \r$.}
\end{equation}
One then defines a Lipschitz continuous bump function $\xi \in C(M)$
around any fixed $x_0 \in M$ by
\begin{equation} \label{xi}
\xi(x) := \eta\left( \frac{d_{\wt{g}}(x,x_0) - (\g + \d/2)}{\d/4} \right), 
\end{equation}
where $d_{\wt{g}}$ is the distance function with respect to the metric $\wt{g}$. 
Note that $d_{\wt{g}}(\cdot,x_0)$ is Lipschitz continuous but need not be smooth 
everywhere, and hence $\xi$ need not be smooth everywhere.
By construction, $\xi$ has the following properties
\begin{equation}
\begin{array}{cl}
  \xi(x) = 1, & x \in  B(x_0,\g+\d/2), \\
  \xi(x) = 0, & x \in M \backslash B(x_0,\g+ 3 \d/4), \\
 \end{array}
\end{equation}
which is illustrated in Figure \ref{fig:CutOff2}. \medskip

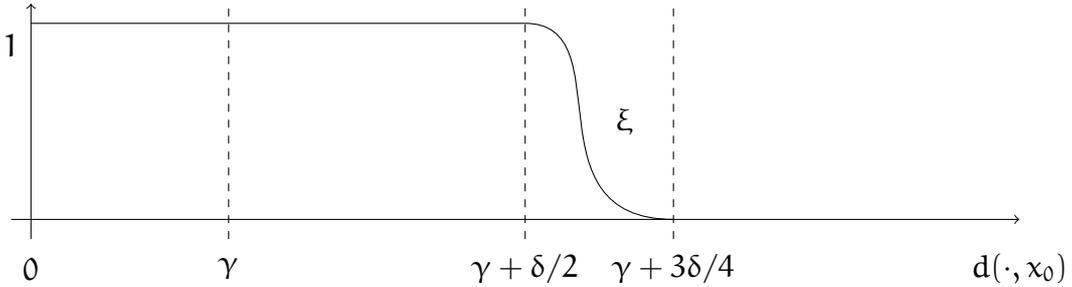
\begin{figure}[h]
	\begin{center}
		
		\begin{tikzpicture}[scale=1.3]
		\draw[->] (1.8,0) -- (12,0);
		\draw[->] (2,-0.2) -- (2,2.2);
		
		\draw (1.8,2) node[anchor=north] {$1$};
		
		\draw (2,2) -- (7,2);
		\draw (7,2) .. controls (8,2) and (7,0) .. (8.5,0);
		\draw[dashed] (4,-0.2) -- (4,2.2);
		\draw[dashed] (7,-0.2) -- (7,2.2);
		\draw[dashed] (8.5,-0.2) -- (8.5,2.2);
		
		\draw (8,1) node {$\xi$};
		
		\draw (2,-0.5) node {$0$};
		\draw (4,-0.5) node {$\g$};
		\draw (7,-0.5) node {$\g+\d/2$};
		\draw (8.5,-0.5) node {$\g+3\d/4$};
		\draw (12,-0.5) node {$d(\cdot,x_0)$};
		
		\end{tikzpicture}
		
		\caption{The bump function $\xi$.}
		\label{fig:CutOff2}
	\end{center}
\end{figure}

\noindent Below in \S \ref{curvature-section}, starting with \eqref{xi-definition}, we provide a careful argument 
differentiating between the case that $\xi$ is 
smooth in a neighborhood of $x$ and the case that $\xi$ is not. The latter case is studied after
\eqref{dg-equation} using a trick of Calabi. In case of smoothness, we have by \eqref{eta properties} 
control on derivatives of $\xi$
\begin{equation} \label{prop cutoff 1}
|\wt{\n} \xi|^2 (x) \le \frac{16^2}{\d^2} \xi(x), \quad x \in M.
\end{equation}

\noindent Shi also proves an estimate
\begin{equation}\label{prop cutoff 2}
\wt{\n} \wt{\n} \xi(x) \ge - \, c_0(\g,\d,k_0)\wt{g}(x), \quad x \in M,
\end{equation}
where $c_0(\g,\d,k_0) > 0$ is a constant only depending on $\g,\d$ and $k_0$. 
\medskip

\noindent The auxiliary bump function $\xi$ is used to define
\[ F(x,t) := \xi(x) \psi(x,t), \quad (x,t) \in B(x_0, \g+\d) \times [0,T]. \]
By construction, it has the properties
\begin{equation}
\begin{split}
&F(x,0) = 0, \quad x \in B(x_0, \g+\d), \\
&F(x,t) = 0, \quad (x,t) \in M \backslash B(x_0, \g+ 3\d/4) \times [0,T], 
\end{split}
\end{equation}
In particular, $F$ attains its maximum on $B(x_0,\g+ 3\d/4) \times [0,T]$, i.e.
there exists $(x_0,t_0) \in B(x_0,\g+ 3\d/4) \times [0,T]$ such that
\[ 
F(x_0,t_0) = \max \, \{ \, F(x,t) \mid (x,t) \in B(x_0,\g+ \d) \times [0,T] \}. 
\]
Using the evolution inequality \eqref{psi equation}, especially the negative quadratic term 
$(- \frac{1}{16} \psi^2)$, as well as the properties \eqref{prop cutoff 1} and \eqref{prop cutoff 2} 
of the cutoff function $\xi$, Shi concludes by maximum principle arguments that
\[ 
F(x_0,t_0) \le c(n,\g,\d,T,\wt{g}), 
\]
where $c(n,\g,\d,T,\wt{g})>0$ is a constant only depending on $n,\g,\d,T,\wt{g}$. Thus
\begin{equation}\label{F-estimate}
\xi(x) \psi(x,t) = F(x,t) \le F(x_0,t_0) \le c(n,\g,\d,T,\wt{g}),  
\end{equation}
for any $(x,t) \in B(x_0, \g+\d) \times [0,T]$. Since
$\xi \equiv 1$ on $B(x_0,\g+\d/2)$, we conclude
\begin{equation}\label{psi estimate}
|\wt{\n}g|^2 \varphi(x,t) = \psi(x,t) \le c(n,\g,\d,T,\wt{g}), 
\end{equation}
for any $(x,t) \in B(x_0, \g+ \d/2) \times [0,T]$. Finally, 
since by definition $\varphi(x,t) \ge a$, the statement follows from
\[ |\wt{\n}g|^2(x,t) \le \frac{1}{a}c(n,\g,\d,T,\wt{g}), \quad (x,t) \in B(x_0, \g+\d/2) \times [0,T]. \]
\end{proof}

\begin{lem}[\cite{shi}, Lemma 4.2]\label{Shi Lemma 4.2}
Under the same assumptions as in \cref{Shi Lemma 4.1}, there exists 
a constant $c(n,m,\g,\d,T,\wt{g}) > 0$ for any $m \ge 0$, depending only on $n,m,\g,\d,T$ and $\wt{g}$, such that
\begin{equation}\label{}
|\wt{\n}^m g(x,t)|^2 \le c(n,m,\g,\d,T,\wt{g})
\end{equation}
for all $(x,t) \in B(x_0,\g+\frac{\d}{m+1}) \times [0,T]$.
\end{lem}

\begin{proof}[Proof outline] 
\cref{Shi Lemma 4.2} is proven by induction. Assuming that the statement 
holds for any integer $0 \leq m_0 < m$, Shi defines the function (cf. \eqref{psi})
\[ \Psi(x,t) = (a_0 + |\wt{\n}^{m-1} g(x,t)|^2) |\wt{\n}^m g(x,t)|^2 \]
and proves that, if $a_0 > 0$, depending only on $m,n,\g,\d,T, \wt{g}$, is chosen appropriately,
then $\Psi$ satisfies a differential inequality of the form (cf. \eqref{psi equation})
\[ 
\frac{\p \Psi}{\p t} \le g^{\a\b} \wt{\n}_\a \wt{\n}_\b \Psi - c_1 \Psi^2 + c_0, 
\]
on $B(x_0,\g + \d/m) \times [0,T]$, where $c_0, c_1 > 0$ only depend on $m,n,\g,\d,T$ and $\wt{g}$. 
Then by the same steps as in the proof of \cref{Shi Lemma 4.1}, Shi obtains (cf. \eqref{psi estimate})
\[ \Psi(x,t) \le c_2(m,n,\g,\d,T,\wt{g}), \quad 
\textup{for} \ (x,t) \in B(\ov{U},\d/(m+1)) \times [0,T]. \]
Hence, we conclude for all $(x,t) \in B(x_0,\d/(m+1)) \times [0,T]$
\[ 
|\wt{\n}^m g(x,t)|^2 \le \frac{1}{a_0} \Psi(x,t) \le \frac{1}{a_0} c_2(m,n,\g,\d,T,\wt{g}),
\]
which finishes the proof.
\end{proof}

Now Shi completes the proof of \cref{Shi Theorems} as follows. Shi takes an exhaustion of the manifold 
$M$ by relatively compact domains $D_k \subset M$, $k \in \N_0$, with $(n-1)$-dimensional, smooth, compact boundary $\p D_k$, 
such that $B(x_0,k) \subset D_k$, for some fixed point $x_0 \in M$. By \cref{Shi Theorem 3.2} and \cref{Shi apriori}, 
there exists $T(n,k_0)>0$ depending only on $n$ and $k_0$ such that
the system (cf. \eqref{Dirichlet Ricci})
\begin{equation} \label{Dirichlet Ricci-k}
\begin{split}
\frac{\p}{\p t} g_{ij}(x,t) = (-2 \Ric_{ij} + \n_i V_j + \n_j V_i)(x,t)&, \quad  (x,t) \in D_k \times [0,T],  \\
g(x,t) = \wt{g}(x)&, \quad (x,t) \in \partial D_k \times [0,T], \\
g(x,0) = \wt{g}(x)&, \quad x \in D_k.
\end{split}
\end{equation}
has a unique smooth solution $g(k,x,t)$ on $D_k \times [0,T(n,k_0)]$ satisfying 
\begin{equation}
(1-\ve(n)) \wt{g}(x) \le g(k,x,t) \le (1+\ve(n)) \wt{g}(x)
\end{equation}
for all $(x,t) \in D_k \times [0,T(n,k_0)]$. Here, $\ve(n)>0$
is a sufficiently small constant, depending only on $n$, introduced in Lemma \ref{Shi Lemma 4.1}. 
Now, for any $k \ge 2$, the solution $g(k,x,t)$ is defined on $B(x_0,1)$. By \cref{Shi Lemma 4.2}, we have for all $m \in \N_0$
\begin{equation}\label{}
|\wt{\n}^m g(k,x,t)|^2 \le c(n,m,q,T(n,k_0),\wt{g})
\end{equation}
for all $(x,t) \in B(x_0,1) \times [0,T(n,k_0)]$ and all $k \ge 2$. Hence by 
Arzel\`a-Ascoli there exists a subsequence $(g(k_\ell,x,t))_{\ell \in \N_0}$, which converges on 
$B(x_0,1) \times [0,T(n,k_0)]$ in the $C^\infty$ topology to a family of smooth metrics $g(x,t)$.
\medskip

By the same argument a subsequence of this subsequence converges on $B(x_0,2) \times [0,T(n,k_0)]$.
We iterate this argument and consider the diagonal sequence.
Then, for every fixed $q \in \N$, the diagonal sequence converges to $g(x,t)$ on $B(x_0,q) \times [0,T(n,k_0)]$, 
and thus converges smoothly locally uniformly to $g(x,t)$. Thus $g(x,t)$ solves \eqref{Ricci de Turck system on M incomplete}. 
The estimate \eqref{unif equiv delta on M complete} follows by restricting the solutions $g(k,x,t)$ to $0 \le t \le T(n,k_0,\d)$, 
where $T(n,k_0,\d)$ is from \cref{Shi apriori}. 

\section{A priori estimates of $\nabla g$ along the flow}\label{curvature-section}

In this section we establish quantitative estimates for the first derivatives of the metric under Ricci de Turck flow on singular manifolds. We assume bounded curvature at time $t=0$ and that the metrics $g(t)$ are uniformly equivalent and sufficiently close to the initial metric $\wt{g}$. As a byproduct we also obtain an estimate on the de Turck vector field $V$. We continue in the setting of an $n$-dimensional, smooth and possibly incomplete Riemannian manifold 
$(M,\wt{g})$ and prove an analogue of Lemma \ref{Shi Lemma 4.1}. 

\begin{lem}\label{est nab g}
	Consider $x_0 \in M$ and fix any\footnote{Below, in Corollary \ref{est nab g cor} we will set 
	$\gamma = \delta>0$ sufficiently small.} finite $\g, \d, T >0$ with $\d \le 1$. 
	Let $g(x,t)$ be a smooth solution of the initial value problem
	\begin{equation*}
	\begin{split}
	\frac{\p}{\p t} g_{ij}(x,t) = (-2 \Ric_{ij} + \n_i V_j + \n_j V_i)(x,t)&, \quad  (x,t) \in \, B(x_0,\g+\d) \times [0,T],  \\
	g(x,0) = \wt{g}(x)&, \quad x \in \, B(x_0,\g+\d),
	\end{split}
	\end{equation*}
	where $V^i = g^{jk}(\G^i_{jk} - \wt{\G}^i_{jk})$ is the de Turck vector field. 
	We assume that $B(x_0,\g+\d - r)$ is relatively compact in $M$ for all $r > 0$. Furthermore, we assume that
	for all $(x,t) \in B(x_0,\g+\d) \times [0,T]$ we have the inequalities
	\begin{equation} \label{equiv met}
	(1-\ve(n)) \wt{g}(x) \le g(x,t) \le (1+\ve(n)) \wt{g}(x) 
	\end{equation}
	for $\ve(n) > 0$ sufficiently small, only depending on $n$. Also assume that 
	\[ |\wt{Rm}|^2 \le k_0 \]
	for some constant $k_0 > 0$. Then there exist constants $c(n), c(n,k_0) > 0$, 
	only depending on the arguments in brackets, such that
	for all $(x,t) \in B(x_0,\g+\frac{\d}{2}) \times [0,T]$
	\begin{align}\label{nab g est lemma}
	\begin{split}
	|\wt{\n}g| (x,t) & \le \frac{c(n,k_0)}{\d} + c(n)c_1, \quad \textup{where} \ \ 
	c_1 := \sup_{x \in B(x_0,\g+3\d/4)} |\wt{\n} \wt{\Rm}|(x).
	\end{split}
	\end{align}
\end{lem}
\begin{rem}
The restriction $\d \le 1$ is for technical reasons to achieve a simpler expression for the right-hand side of \eqref{nab g est lemma}. For our purposes this is sufficient as we are aiming at estimates on an incomplete manifold when we get closer and closer to the singularity.	
Also note that the estimates \eqref{nab g est lemma} are independent of $\gamma$, and only depend on the difference of radia of the smaller
ball $B(x_0,\g+\frac{\d}{2})$ and the larger ball $B(x_0,\g+\d)$.
\end{rem}

We will prove the lemma below and first note its consequence $-$
estimates on the first derivatives of the metric for Ricci de Turck flow. More specifically, 
assuming additionally that $|\wt{\n} \wt{\Rm}| = \mathcal{O}(\rho^{-1})$, where $\rho > 0$ is the distance to the singularity, 
a natural condition in case $|\wt{\Rm}|$ is bounded, we obtain that $|\wt{\n} g| = \mathcal{O}(\rho^{-1})$ and 
$|V| = \mathcal{O}(\rho^{-1})$ uniformly in $t \in [0,T]$.
\begin{cor}\label{est nab g cor}
Let $(M,\wt{g})$ be a (possibly incomplete) smooth Riemannian manifold of dimension $n$. 
Fix $0 < T < \infty$ and let $g(x,t)$ be a smooth solution of 
\begin{equation*}
\begin{split}
\frac{\p}{\p t} g_{ij}(x,t) = (-2 \Ric_{ij} + \n_i V_j + \n_j V_i)(x,t)&, \quad  (x,t) \in \, M \times [0,T],  \\
g(x,0) = \wt{g}(x)&, \quad x \in \, M,
\end{split}
\end{equation*}
where $V$ is the de Turck vector field as above. Assume that for all $(x,t) \in M \times [0,T]$
\begin{equation*}
(1-\ve(n)) \wt{g}(x) \le g(x,t) \le (1+\ve(n)) \wt{g}(x) 
\end{equation*}
for $\ve(n) > 0$ sufficiently small, only depending on $n$, and also assume that there exist constants $k_0, C > 0$, such that 
\[ |\wt{Rm}|^2 \le k_0 \] 
and that for all $x \in M$, $0 < \rho \le 1$
\[ |\wt{\n}\wt{\Rm}|(x) \le \frac{C}{\rho}  \] 
whenever $B(x,\rho-r)$ is relatively compact for all $r > 0$. 
Then there exists $C' > 0$ such that for all $x \in M$, $t \in [0,T]$, $0 < \rho \le 1$
\[ |\wt{\n}g|(x,t) \le \frac{C'}{\rho}, \quad  |V|(x,t) \le \frac{C'}{\rho} \]
whenever $B(x,\rho-r)$ is relatively compact for all $r > 0$.
\end{cor}	

\begin{rem}
The (technical) condition $B(x,\rho-r)$ is relatively compact in $M$ for all $r > 0$ is a way to express the distance of a point $x \in M$ to the singular strata of $M$ intrinsically. It means that this distance is larger or equal to $\rho$. 	
\end{rem}

\begin{proof}[Proof of \cref{est nab g cor}]
	Consider $x_0 \in M$ and $\rho \le 1$ such that $B(x_0,\rho-r)$ is relatively compact in $M$ 
	for all $r > 0$. Then by \cref{est nab g} (choosing $\g,\d$ in \cref{est nab g} as equal to $\rho/2$) we obtain
	\[ |\wt{\n}g|(x_0,t) \le \frac{c(n,k_0)}{\rho} + c(n) c_1, \]
	where the constant $c_1$ can be estimated as follows
	\[ c_1 = \sup_{x \in B(x_0,7\d/8)} |\wt{\n}\wt{\Rm}|(x) \le \frac{8C}{\rho},\]
	since for all $x \in B(x_0,7\rho/8)$ we have that $B(x,\rho/8-r)$ is relatively 
	compact for all $r > 0$. This proves the estimate for $ |\wt{\n}g|$.
	The estimate of the de Turck vector field $V$ follows from this and
	\[ V = g^{-1} * \wt{\n}g, \]
	see \cite[p. 266, formula (32)]{shi}.
	\end{proof}

\noindent We can now proceed with proof of \cref{est nab g}.

\begin{proof}[Proof of \cref{est nab g}]
Our strategy is a careful analysis of the proof of \cite[Lemma 4.1]{shi},
which is written out here in Lemma \ref{Shi Lemma 4.1}, while making the 
dependencies of various constants explicit. For the convenience of the reader 
and to keep our argument here self-contained, we repeat the steps from \cite[Lemma 4.1]{shi} here. \medskip
 
\noindent In the following, $c(n)$ and $c(n,k_0)$ denote constants only depending on 
$n$ and $n, k_0$, respectively. The constants may vary from estimate to estimate.\medskip 

\noindent As in \cite[Proof of Lemma 4.1, p.247 (5)]{shi} we have 
\begin{align}
\begin{split}
\frac{\p}{\p t} |\wt{\n}g|^2 = & g^{\a\b} \wt{\n}_\a \wt{\n}_\b |\wt{\n} g|^2 - 2 g^{\a\b} \wt{\n}_\a \wt{\n} g \cdot \wt{\n}_\b \wt{\n} g \\
& + \wt{\Rm} * g^{-2} * g * \wt{\n} g * \wt{\n} g + g^{-1} * g * \wt{\n} \wt{\Rm} * \wt{\n} g \\
& + g^{-2} * \wt{\n} g * \wt{\n} g * \wt{\n} \wt{\n} g + g^{-3} * \wt{\n} g * \wt{\n} g * \wt{\n} g * \wt{\n} g.
\end{split}
\end{align}
Here the product $A*B$ of two tensors $A$ and $B$ denotes a linear combination of terms which are obtained as follows: Starting from the tensor product $A \otimes B$, perform an arbitrary number of the following operations: taking contractions, raising, lowering or permuting indices. The important consequence in our case here is that it will always be possible to estimate
\[ 
|A * B| \le c(n) |A| \cdot |B|,
\]
where $c(n)$ depends on the specific form of the product. 
Since by assumption, the closure $\ov{B(x_0,\g + \frac{3}{4}\d)} \subset M$ is compact, we conclude
($c_1$ is defined in \eqref{nab g est lemma})
\begin{align}\label{nab rm bound}
|\wt{\n}\wt{\Rm}| \le c_1 \ \text{ on } B(x_0,\g + \frac{3}{4}\d).
\end{align}
Furthermore, by \eqref{equiv met} we have
\begin{align}\label{equiv met 2}
\frac{1}{2}\wt{g}(x) \le g(x,t) \le 2 \wt{g}(x) \quad \text{ on } B(x_0,\g+\d).
\end{align}
Hence 
\begin{align}
\begin{split}
\wt{\Rm} * g^{-2} * g * \wt{\n} g * \wt{\n} g \le c(n,k_0) |\wt{\n} g|^2, \\
g^{-1} * g * \wt{\n} \wt{\Rm} * \wt{\n} g \le c(n) c_1 |\wt{\n} g|
\end{split}
\end{align}
on $B(x_0,\g+3\d/4) \times [0,T]$. Also, whenever we use the bound \eqref{nab rm bound} 
on $\wt{\n} \wt{\Rm}$ it is understood that the estimate, which follows, holds on $B(x_0,\g+3\d/4) \times [0,T]$. 
As in \cite[Proof of Lemma 4.1, p.247 (9)]{shi} we have
\begin{align}
\begin{split}
g^{-2} * \wt{g} * \wt{\n} g * \wt{\n} \wt{\n} g \le 72n^5 |\wt{\n} g|^2 |\wt{\n}^2 g|, \\
g^{-3} * \wt{\n} g * \wt{\n} g * \wt{\n} g * \wt{\n} g \le 160 n^6 |\wt{\n} g|^4.
\end{split}
\end{align}
This gives
\begin{align}\label{inbetween 1}
\begin{split}
\frac{\p}{\p t} |\wt{\n} g|^2 \le & g^{\a\b} \wt{\n}_\a \wt{\n}_\b |\wt{\n} g|^2 - |\wt{\n}^2 g|^2 + c(n,k_0) |\wt{\n} g|^2 + c(n) c_1 |\wt{\n}g| \\ 
& + 72n^5 |\wt{\n} g|^2 |\wt{\n}^2 g| + 160 n^6 |\wt{\n} g|^4. 
\end{split}
\end{align}
Estimating as in \cite[Proof of Lemma 4.1, p.247]{shi}
\begin{equation}\label{inbetween 2}
\begin{split}
72n^5 |\wt{\n} g|^2 |\wt{\n}^2 g| + 160 n^6 |\wt{\n} g|^4 &\le \frac{1}{2} |\wt{\n}^2 g|^2 + 3200n^{10} |\wt{\n} g|^4,
\\ c(n) c_1 |\wt{\n}g| &\leq \frac{(c(n) c_1)^2}{2} + \frac{|\wt{\n}g|^2}{2},
\end{split}
\end{equation}
we obtain from \eqref{inbetween 1} after an appropriate change of constants $c(n,k_0)$ and $c(n)$
\begin{equation}\label{norm nab g inequ}
\begin{split}
\frac{\p}{\p t} |\wt{\n} g|^2  &\le  g^{\a\b} \wt{\n}_\a \wt{\n}_\b |\wt{\n} g|^2 - \frac{1}{2}|\wt{\n}^2 g|^2 + 3200n^{10} |\wt{\n} g|^4 
\\ &+ c(n,k_0) |\wt{\n} g|^2 + c(n) c_1^2.
\end{split}
\end{equation}
As in \cite[Proof of Lemma 4.1, p.248]{shi}, we fix a small constant $\ve \equiv \ve(n) := (256000n^{10})^{-1}$, 
such that the inequality \eqref{equiv met} now reads as
\begin{equation}\label{lambda-eigenvalue-estimates}
1-\ve(n) \le \l_k(x,t) \le 1+\ve(n),
\end{equation}
for any $k = 1, 2, \dots, n$, where $\l_k(x,t)$ refers to the eigenvalues of $g(x,t)$ with respect to $\wt{g}(x)$. 
Sometimes we use a rougher estimate $\frac{1}{2} \le \l_k(x,t) \le 2$ instead. We also set 
\begin{equation}
m := 25600n^{10}, \qquad a := 6400n^{10}
\end{equation}
and define (we simplify notation by writing $\l_k \equiv \l_k(x,t)$)
\begin{equation}
\varphi(x,t) := a + \sum_{k=1}^n \l_k^m, \qquad (x,t) \in B(x_0,\g+\d) \times [0,T].
\end{equation}
Following \cite[Proof of Lemma 4.1, p.248 (16)]{shi} we obtain
\begin{align}\label{dphi-dt-estimate}
\begin{split}
\frac{\p \varphi}{\p t} = & m \l_k^{m-1} g^{\a\b} \wt{\n}_\a \wt{\n}_\b g_{kk} \\
& + m \l_k^{m-1} * (\wt{\Rm} * g^{-1} * g + g^{-2} * \wt{\n} g * \wt{\n} g).
\end{split}
\end{align}
We now proceed as in Lemma \ref{Shi Lemma 4.1} along the following steps.
\medskip

\begin{enumerate}
\item Step 1: Derive an evolution inequality for $\psi := \varphi \cdot |\wt{\n} g|^2$ as in \eqref{psi equation}. 
\smallskip

\item Step 2: Estimate $\wt{\n} \wt{\n} \xi$ from below as in \eqref{prop cutoff 2}.
\smallskip

\item Step 3: Estimate $\xi \psi$ from above as in \eqref{F-estimate} and conclude the proof.
\medskip

\end{enumerate}

\noindent \textbf{Step 1: Derive an evolution inequality for $\psi := \varphi \cdot |\wt{\n} g|^2$ as in \eqref{psi equation}. }
\medskip

\noindent We estimate the individual terms on the right hand side of \eqref{dphi-dt-estimate}
\begin{align}
\begin{split}
m \l_k^{m-1} * \wt{\Rm} * g^{-1} * g &\le c(n,k_0), \\
m \l_k^{m-1} * g^{-2} * \wt{\n} g * \wt{\n} g &\le 10n^3 m (1+\ve)^{m-1} |\wt{\n} g|^2.
\end{split}
\end{align}
As in \cite[Proof of Lemma 4.1, p.248]{shi} we have
\begin{align}
\begin{split}
g^{\a\b}\wt{\n}_\a \wt{\n}_\b \varphi & = m \l_k^{m-1} g^{\a\b}\wt{\n}_\a \wt{\n}_\b g_{kk}\\
& \hspace{0.5cm}+ m (\l_i^{m-2} + \l_i^{m-3} \l_j + \cdots + \l_j^{m-2}) \cdot g^{\a\b} \wt{\n}_\a g \cdot \wt{\n}_\b g \\
& \ge m \l_k^{m-1} g^{\a\b}\wt{\n}_\a \wt{\n}_\b g_{kk} + \frac{m(m-1)}{2} (1-\ve)^{m-2}|\wt{\n} g|^2.
\end{split}
\end{align}
This yields
\begin{align}\label{dphi-dt-estimate1}
\begin{split}
\frac{\p \varphi}{\p t} \le & g^{\a\b} \wt{\n}_\a \wt{\n}_\b \varphi - \frac{m(m-1)}{2} (1-\ve)^{m-2}|\wt{\n} g|^2 \\
& + c(n,k_0) + 10n^3 m (1+\ve)^{m-1} |\wt{\n} g|^2.
\end{split}
\end{align}
As in \cite[p.249 (20),(21),(22)]{shi}, we easily check
\begin{align}\label{mn-estimates-algebra}
\begin{split}
&10n^3 m (1+\ve)^{m-1} \le \frac{m^2}{16}, \\
&\frac{m(m-1)}{2} (1-\ve)^{m-2} \ge \frac{m^2}{4}(1-\ve)^{m-2} \ge \frac{3}{16}m^2,
\end{split}
\end{align}
such that \eqref{dphi-dt-estimate1} reduces to
\begin{align}\label{phi inequ}
\frac{\p \varphi}{\p t} \le & g^{\a\b} \wt{\n}_\a \wt{\n}_\b \varphi + c(n,k_0) - \frac{m^2}{8}  |\wt{\n} g|^2.
\end{align}
From \eqref{norm nab g inequ} and \eqref{phi inequ} it follows that
\begin{align}\label{dphi-dt-estimate2} 
\begin{split}
\frac{\p}{\p t}(\varphi \cdot |\wt{\n} g|^2) \le & g^{\a\b} \wt{\n}_\a \wt{\n}_\b 
(\varphi \cdot |\wt{\n} g|^2) -2g^{\a\b} \wt{\n}_\a \varphi \wt{\n}_\b |\wt{\n} g|^2 - \frac{\varphi}{2} |\wt{\n}^2 g|^2 \\
& + 3200n^{10} \varphi |\wt{\n} g|^4 + c(n,k_0) \varphi |\wt{\n}g|^2 + c(n)c_1^2 \varphi \\
& + c(n,k_0) |\wt{\n} g|^2 - \frac{m^2}{8} |\wt{\n} g|^4.
\end{split}
\end{align}
We estimate some of the terms on the right hand side of \eqref{dphi-dt-estimate2}. 
As in \cite[Proof of Lemma 4.1, p.249 (26), p.250 (28)]{shi} we find for the fourth term on the
right hand side of \eqref{dphi-dt-estimate2} 
\begin{align}
\begin{split}
3200n^{10} \varphi |\wt{\n} g|^4 \le 3200n^{10} (a + n(1+\ve)^m) |\wt{\n} g|^4 \le \frac{m^2}{16} |\wt{\n} g|^4.
\end{split}
\end{align}
The second term on the right hand side of \eqref{dphi-dt-estimate2} is estimated as follows.
\begin{align}
\begin{split}
-2 g^{\a\b} \wt{\n}_\a \varphi \wt{\n}_\b |\wt{\n}g|^2 & = -2 g^{\a\b} \wt{\n}_\a \left(\sum_{k=1}^n \l_k^m\right) \cdot \wt{\n}_\b |\wt{\n}g|^2 \\
& = -4 g^{\a\b} \cdot \left( m \, \l_k^{m-1} \cdot \wt{\n}_\a \l_k \right) \cdot \wt{\n}_\b |\wt{\n}g|^2 \\
& \le 8mn^5 (1+\ve)^{m-1} |\wt{\n}g|^2 |\wt{\n}^2 g| \\
& \le \sqrt{\phi} \, |\wt{\n}^2 g| \cdot \left( \frac{16mn^5 |\wt{\n}g|^2}{\sqrt{\phi}} \right) \\
& \le \frac{\varphi}{2} |\wt{\n}^2 g|^2 + \frac{128m^2 n^{10}}{\varphi} |\wt{\n} g|^4.
\end{split}
\end{align}
Plugging these estimates back into \eqref{dphi-dt-estimate2} yields
\begin{align}\label{dphi-dt-estimate3}
\begin{split}
\frac{\p}{\p t}(\varphi \cdot |\wt{\n} g|^2) \le & g^{\a\b} \wt{\n}_\a \wt{\n}_\b (\varphi \cdot |\wt{\n} g|^2) + \frac{128m^2 n^{10}}{\varphi} |\wt{\n} g|^4 - \frac{m^2}{16} |\wt{\n} g|^4\\
& + c(n,k_0) \varphi |\wt{\n}g|^2 + c(n)c_1^2.
\end{split}   
\end{align}
Since $\varphi(x,t) \ge a$, with $a = 6400n^{10}$, we have 
\[ \frac{128m^2 n^{10}}{\varphi} \le \frac{m^2}{32}, \]
such that \eqref{dphi-dt-estimate3} reduces to
\begin{align}\label{dphi-dt-estimate4}
\begin{split}
\frac{\p}{\p t}(\varphi \cdot |\wt{\n} g|^2) \le & g^{\a\b} \wt{\n}_\a \wt{\n}_\b (\varphi \cdot |\wt{\n} g|^2) - \frac{m^2}{32} |\wt{\n} g|^4\\
& + c(n,k_0) \varphi |\wt{\n}g|^2 + c(n)c_1^2.
\end{split}   
\end{align}
Using \eqref{lambda-eigenvalue-estimates} and the first estimate of \eqref{mn-estimates-algebra}
in the second inequality, we find
\begin{align}
\frac{m^2}{32} |\wt{\n}g|^4 \equiv  \frac{m^2}{32 \varphi} |\wt{\n}g|^4 \varphi 
\ge \frac{m^2}{32 (a+n(1+\ve)^m)^2} |\wt{\n}g|^4 \varphi^2 \ge \frac{1}{8} |\wt{\n}g|^4 \varphi^2.
\end{align}
Thus we obtain from \eqref{dphi-dt-estimate4}, using the inequality $ab \leq \frac{1}{2} a^2+ \frac{1}{2} b^2$ 
and adapting the constant $c(n,k_0)>0$ accordingly in the last estimate
\begin{align}
\begin{split}
\frac{\p}{\p t}(\varphi \cdot |\wt{\n} g|^2) & \le g^{\a\b} \wt{\n}_\a \wt{\n}_\b (\varphi \cdot |\wt{\n} g|^2) - \frac{1}{8} |\wt{\n} g|^4 \varphi^2\\
& \hspace{0.5cm} + c(n,k_0) \varphi |\wt{\n}g|^2 + c(n)c_1^2 \\
& \le g^{\a\b} \wt{\n}_\a \wt{\n}_\b (\varphi \cdot |\wt{\n} g|^2) - \frac{1}{16} |\wt{\n} g|^4 \varphi^2\\
& \hspace{0.5cm} + c(n,k_0) + c(n)c_1^2.
\end{split}   
\end{align}
Defining $\psi(x,t) := (\varphi \cdot |\wt{\n} g|^2)(x,t)$ this inequality reads
\begin{align}\label{psi inequ}
\begin{split}
\frac{\p \psi}{\p t} & \le g^{\a\b} \wt{\n}_\a \wt{\n}_\b \psi  - \frac{1}{16} \psi^2 + c(n,k_0) + c(n)c_1^2.
\end{split}   
\end{align} \ \medskip

\noindent \textbf{Step 2: Estimate $\wt{\n} \wt{\n} \xi$ from below as in \eqref{prop cutoff 2}. }
\medskip

\noindent Next, as in \cite[Proof of Lemma 4.1, p.251 (36),(37)]{shi} we take a cutoff function $\eta \in C^\infty(\r)$ 
as in \eqref{eta properties}, illustrated in Figure \ref{fig:CutOff1}. 
Then we define the cutoff function $\xi \in C(M)$
\begin{align}\label{xi-definition}
\begin{split}
\xi(x) = \eta \left(\frac{d_{\wt{g}}(x,x_0) - (\g+\d/2)}{\d/4}\right),
\end{split}
\end{align}
where $d_{\wt{g}}$ is the distance function with respect to the metric $\wt{g}$. Note that $d_{\wt{g}}(\cdot,x_0)$ is Lipschitz continuous but need not be smooth everywhere, and hence $\xi$ need not be smooth everywhere. From the properties of $\eta$ we have 
\begin{align}
\begin{split}
& \xi(x) = 1, \qquad x \in B(x_0,\g+\d/2), \\
& \xi(x) = 0, \qquad x \in M \backslash B(x_0,\g+3\d/4), \\
& 0 \le \xi(x) \le 1, \qquad x \in M.
\end{split}
\end{align} 
If $d_{\wt{g}}(\cdot,x_0)$ is smooth in a neighborhood of a point $x$, then we also have
\begin{align}
\wt{\n}_\b \xi(x) = \frac{4}{\d} \eta'\left(\frac{d_{\wt{g}}(x,x_0) - (\g+\d/2)}{\d/4}\right) \wt{\n}_\b d_{\wt{g}}(x,x_0) 
\end{align}
\begin{align}\label{step3-estimate1}
\begin{split}
\wt{\n}_\a \wt{\n}_\b \xi(x) = & \frac{4}{\d} \eta'\left(\frac{d_{\wt{g}}(x,x_0) - (\g+\d/2)}{\d/4}\right) \wt{\n}_\a \wt{\n}_\b d_{\wt{g}}(x,x_0) \\
& + \frac{16}{\d^2} \eta''\left(\frac{d_{\wt{g}}(x,x_0) - (\g+\d/2)}{\d/4}\right) \wt{\n}_\a d_{\wt{g}}(x,x_0) \wt{\n}_\b d_{\wt{g}}(x,x_0). 
\end{split}
\end{align}
Since $|\wt{\n} d_{\wt{g}}(x,x_0)| = 1$, it follows using $|\eta'|^2 \le 16 \eta$ that
\begin{align} \label{inbetween 3}
|\wt{\n} \xi(x)|^2 \le \frac{16}{\d^2} (\eta')^2\left(\frac{d_{\wt{g}}(x,x_0) - (\g+\d/2)}{\d/4}\right) \le \frac{256}{\d^2} \xi(x).
\end{align}
Furthermore, note that
\begin{align}
\wt{\n}_\a d_{\wt{g}}(x,x_0) \wt{\n}_\b d_{\wt{g}}(x,x_0) \le \wt{g}_{\a\b}(x),
\end{align}
such that, using $|\eta''| \le 8$, we can estimate from below
\begin{align}\label{step3-estimate2}
\frac{16}{\d^2} \eta''\left(\frac{d_{\wt{g}}(x,x_0) - (\g+\d/2)}{\d/4}\right) \wt{\n}_\a d_{\wt{g}}(x,x_0) \wt{\n}_\b d_{\wt{g}}(x,x_0) \ge -\frac{128}{\d^2} \wt{g}_{\a\b}(x).
\end{align}
By assumption, $|\wt{\Rm}|^2 \le k_0$ and thus the sectional curvature is in particular bounded from below
$\sec \ge -\sqrt{k_0}$. From the Hessian comparison theorem, applied in a relatively compact ball, we conclude
\begin{align}
\wt{\n}_\a \wt{\n}_\b d_{\wt{g}}(x,x_0) \le \sqrt[4]{k_0} \coth \left(\sqrt[4]{k_0} d_{\wt{g}}(x,x_0) \right) \wt{g}_{\a\b}(x).
\end{align}
Using $0 \ge \eta'(s) \ge -4\, \eta^{1/2}(s) \ge -\, 4$ for all $s\in \r$,
it follows that 
\begin{equation}\label{step3-estimate3}
\begin{split}
&\frac{4}{\d}\eta'\left(\frac{d_{\wt{g}}(x,x_0) - (\g+\d/2)}{\d/4}\right) \wt{\n}_\a \wt{\n}_\b d_{\wt{g}}(x,x_0) 
\\ &\ge -\frac{16}{\d} \sqrt[4]{k_0} \coth \left(\sqrt[4]{k_0} d_{\wt{g}}(x,x_0) \right) \wt{g}_{\a\b}(x).
\end{split}
\end{equation}
We now obtain from \eqref{step3-estimate1}, combined with \eqref{step3-estimate2} and \eqref{step3-estimate3}
\begin{align}\label{hess est xi}
\begin{split}
\wt{\n}_\a \wt{\n}_\b \xi(x) \ge -\left( \frac{128}{\d^2} + \frac{16}{\d} \sqrt[4]{k_0} \coth \left(\sqrt[4]{k_0} d_{\wt{g}}(x,x_0) \right) \right) \wt{g}_{\a\b}(x).
\end{split}
\end{align} \ \medskip

\noindent \textbf{Step 3: Estimate $\xi \psi$ from above as in \eqref{F-estimate} and conclude the proof.}
\medskip

\noindent 
Next we simplify notation by writing as in the proof of Lemma \ref{Shi Lemma 4.1}
\[ F(x,t) := \xi(x)\psi(x,t), \qquad (x,t) \in B(x_0,\g+\d) \times [0,T]. \]
Since $|\wt{\n}g|^2(x,0) = 0$, we have
\begin{align}
F(x,0) = 0, \qquad x \in B(x_0,\g+\d).
\end{align}
Since $\xi(x) = 0$ for $x \in B(x_0,\g+\d) \backslash B(x_0,\g+\frac{3}{4}\d)$, it follows that
\begin{align}
F(x,t) = 0, \qquad (x,t) \in B(x_0,\g+\d) \backslash B(x_0,\g+\frac{3}{4}\d) \times [0,T].
\end{align}
Thus there exists a point $(y_0,t_0) \in B(x_0,\g+\frac{3}{4}\d) \times [0,T]$ with $t_0 > 0$ such that
\begin{align}
F(y_0,t_0) = \max \ \left\{F(x,t) \mid (x,t) \in B(x_0,\g+\d) \times [0,T] \right\}
\end{align}
unless $F \equiv 0$ on $B(x_0,\g+\d) \times [0,T]$. \medskip

In the following, as already alluded to in the proof of Lemma \ref{Shi Lemma 4.1},
we distinguish three cases, first case where $\xi \equiv 1$ in a neighborhood of $y_0$, 
second case where $\xi$ is not identically $1$, but smooth in a neighborhood of $y_0$, 
and third case, where $\xi$ is not smooth and a trick needs to be applied. \medskip

\noindent \underline{\textbf{Case 1. $y_0 \in B(x_0,\g+\d/2)$}} \medskip

\noindent Then $\xi \equiv 1$ near $y_0$, such that $F = \psi$ near $(y_0,t_0)$, and we have by \eqref{psi inequ}
\begin{align}
0 \le (\frac{\p}{\p t} - g^{\a\b}\wt{\n}_\a \wt{\n}_\b) \psi(y_0,t_0) \le -\frac{1}{16} \psi^2(y_0,t_0) + c(n,k_0) + c(n)c_1^2, 
\end{align} 
and thus we conclude
\begin{align}
\frac{1}{16}F^2(y_0,t_0) = \frac{1}{16}\psi^2(y_0,t_0) \le c(n,k_0) + c(n)c_1^2.
\end{align}
This estimate is better than the one we will obtain in Case 2. \\

\noindent \underline{\textbf{Case 2. $y_0 \notin B(x_0,\g+\d/2)$ and $y_0$ is not in the cut 
locus\footnote{with respect to the metric $\wt{g}$} of $x_0$}}
\medskip

\noindent Then the distance function $d_{\wt{g}}(\cdot,x_0)$, and hence also $\xi$, is smooth in a neighborhood of 
$y_0$ and it follows that
\begin{align}\label{inequ at max}
\begin{split}
& 0 \le \frac{\p F}{\p t}(y_0,t_0) = \xi(y_0) \frac{\p \psi}{\p t}(y_0,t_0), \\
& 0 = \wt{\n}_\a F(y_0,t_0) = (\xi \wt{\n}_\a \psi + \psi \wt{\n}_\a \xi)(y_0,t_0), \\
& 0 \ge g^{\a\b} \wt{\n}_\a \wt{\n}_\b F(y_0,t_0) = (\xi g^{\a\b} \wt{\n}_\a \wt{\n}_\b \psi + \psi g^{\a\b} \wt{\n}_\a \wt{\n}_\b \xi 
\\ &+ 2 g^{\a\b} \wt{\n}_\a \xi \wt{\n}_\b \psi)(y_0,t_0).
\end{split}
\end{align}
Using \eqref{psi inequ} in the final step, we obtain at the point $(y_0,t_0)$
 \begin{align*}
 0 &\le \left( \frac{\p F}{\p t} - g^{\a\b} \wt{\n}_\a \wt{\n}_\b F\right) (y_0,t_0)
 \\ &\le \xi(y_0) \left( \frac{\p \psi}{\p t} - g^{\a\b} \wt{\n}_\a \wt{\n}_\b \psi \right) (y_0,t_0) 
 \\ &- \left(\psi g^{\a\b} \wt{\n}_\a \wt{\n}_\b \xi + 2 g^{\a\b} \wt{\n}_\a \xi \wt{\n}_\b \psi\right)(y_0,t_0)
 \\ &\le \left( - \frac{1}{16}\xi \psi^2 - \psi g^{\a\b} \wt{\n}_\a \wt{\n}_\b \xi - 2 g^{\a\b} \wt{\n}_\a \xi \wt{\n}_\b \psi \right) (y_0,t_0)
 \\ &+ \xi(y_0) \Bigl(c(n,k_0) + c(n)c_1^2\Bigr).
 \end{align*}
 Thus we conclude at the point $(y_0,t_0)$
 \begin{align}\label{case2-estimate1}
 \begin{split}
 \frac{1}{16}\xi \psi^2  \le -\psi g^{\a\b} \wt{\n}_\a \wt{\n}_\b \xi - 2 g^{\a\b} \wt{\n}_\a \xi \wt{\n}_\b \psi + 
 \xi (c(n,k_0) + c(n)c_1^2).
 \end{split}
 \end{align}
 From the second identity in \eqref{inequ at max} in the first step, and using \eqref{inbetween 3} in the second estimate, 
 we obtain at $(y_0,t_0)$
 \begin{align}\label{case2-estimate2}
 - 2 g^{\a\b} \wt{\n}_\a \xi \wt{\n}_\b \psi  = \frac{2\psi}{\xi} g^{\a\b} \wt{\n}_\a \xi \wt{\n}_\b \xi \le \frac{1024}{\d^2} \psi,
 \end{align}
Furthermore, we obtain from \eqref{hess est xi} at $(y_0,t_0)$
 \begin{align}
 -\psi g^{\a\b} \wt{\n}_\a \wt{\n}_\b \xi \le 2n \left( \frac{128}{\d^2} + \frac{16}{\d} \sqrt[4]{k_0} \coth \left(\sqrt[4]{k_0} \g(y_0,x_0) \right) \right) \psi.
 \end{align}
We estimate the $\coth$-term: Since $y_0 \notin B(x_0,\g+\d/2)$ and $\coth$ is monotonically decreasing on the positive real axis, $\coth(\sqrt[4]{k_0} \g(y_0,x_0)) \le \coth(\sqrt[4]{k_0} \d/2)$. Also, since $z \coth z \le 1 + Cz$ for $z > 0$ and since $\d \le 1$
\begin{align}\label{est coth}
\sqrt[4]{k_0}\coth(\sqrt[4]{k_0} \d/2) = \frac{2}{\d} \sqrt[4]{k_0} \frac{\d}{2}\coth(\sqrt[4]{k_0} \d/2) \le \frac{2}{\d} (1+C \sqrt[4]{k_0} \frac{\d}{2}) \le \frac{c(n,k_0)}{\d}.
\end{align}
Thus we obtain at $(y_0,t_0)$
\begin{align}\label{case2-estimate3}
-\psi g^{\a\b} \wt{\n}_\a \wt{\n}_\b \xi \le \frac{c(n,k_0)}{\d^2} \psi.
\end{align}
Plugging \eqref{case2-estimate2} and \eqref{case2-estimate3} into \eqref{case2-estimate1}, leads to
 \begin{align}
 \begin{split}
 \frac{1}{16}\xi \psi^2 \le \frac{c(n,k_0)}{\d^2} \psi + \xi (c(n,k_0) + c(n)c_1^2).
 \end{split}
 \end{align}
Multiplying this inequality with $\xi$ and using $0 \le \xi \le 1$ we obtain
 \begin{align}
\begin{split}
F(y_0,t_0)^2 \le \frac{c(n,k_0)}{\d^2} F(y_0,t_0) + c(n,k_0) + c(n)c_1^2.
\end{split}
\end{align}  
Thus 
\begin{align}\label{F max est}
\begin{split}
F(y_0,t_0) & \le \frac{c(n,k_0)}{\d^2} + c(n,k_0) + c(n)c_1^2 \\
& \le \frac{c(n,k_0)}{\d^2} + c(n)c_1^2,
\end{split}
\end{align} 
assuming $c(n,k_0) \ge 1$ if necessary. Thus 
\begin{align}
\begin{split}
F(x,t) & \le F(y_0,t_0) \le \frac{c(n,k_0)}{\d^2} + c(n)c_1^2
\end{split}
\end{align} 
for all $(x,t) \in B(x_0,\g+\d) \times [0,T]$. Since
\begin{align*} 
&F(x,t) = \xi(x) \varphi(x,t) |\wt{\n}g|^2 (x,t), \\
&\xi(x) = 1 \text{ for } x \in B(x_0,\g+\d/2), \\
&\varphi(x,t) \ge a = 6400n^{10} \text{ for } (x,t) \in B(x_0,\g+\d) \times [0,T], 
\end{align*}
we obtain
\begin{align}
\begin{split}
|\wt{\n}g|^2 (x,t) & \le \frac{1}{6400n^{10}} \left(\frac{c(n,k_0)}{\d^2} + c(n)c_1^2 \right) \\
& = \frac{c(n,k_0)}{\d^2} + c(n)c_1^2
\end{split}
\end{align}
for all $(x,t) \in B(x_0,\g+\d/2) \times [0,T]$. Hence
\begin{align}\label{dg-equation}
\begin{split}
|\wt{\n}g| (x,t) & \le \frac{c(n,k_0)}{\d} + c(n)c_1
\end{split}
\end{align}
for all $(x,t) \in B(x_0,\g+\d/2) \times [0,T]$. \\

\noindent \underline{\textbf{Case 3. $y_0 \notin B(x_0,\g+\d/2)$ and $y_0$ is in the cut locus of $x_0$}}\medskip

\noindent Then we apply Calabi's trick (see e. g. \cite[p.395]{CLN}). Let $c: [0,d_{\wt{g}}(x_0,y_0)] \to M$ be a minimal geodesic from $x_0$ to $y_0$. Note that since $y_0 \in B(x_0,\g+\frac{3}{4}\d)$ the assumption $B(x_0,\g+\d-r) \subset \subset M$ for all $r > 0$ ensures that such a minimal geodesic exists. Fix $\ve > 0$ sufficiently small and define
\[ \xi_\ve(x) := \eta \left( \frac{d_{\wt{g}}(x,c(\ve)) + \ve - (\g+\d/2)}{\d/4} \right), \quad
F_\ve(x,t) := \xi_\ve(x) \psi(x,t).\]
Since $d_{\wt{g}}(x,x_0) \le d_{\wt{g}}(x,c(\ve)) + \ve$ by the triangle inequality and since $\eta$ is monotonically decreasing, we have 
\[ \xi_\ve(x) \le \xi(x) \]
for all $x \in M$. As $d_{\wt{g}}(y_0,x_0) = d_{\wt{g}}(y_0,c(\ve)) + \ve$, we have
$\xi_\ve(y_0) = \xi(y_0)$. Hence 
\begin{align}
\begin{split}
& F_\ve(x,t) \le F(x,t) \qquad \forall x \in B(x_0,\g+\d) \times [0,T], \\
& F_\ve(y_0,t_0) = F(y_0,t_0),
\end{split}
\end{align}
such that $F_\ve$ has a maximum at $(y_0,t_0)$ as well. The point now is that $\g(\cdot,c(\ve))$ is smooth in a neighborhood of $y_0$. Note that the argument for this in the complete case (see \cite[Proof of Lemma 42, p.284]{pet}) also works in our case since $y_0 \in B(x_0,\g+\frac{3}{4}\d)$ and $B(x_0,\g+\d-r) \subset \subset M$ for all $r > 0$. Thus $F_\ve$ is smooth in a neighborhood of $(y_0,t_0)$ and we can apply the same steps as in Case 2 to $F_\ve$. Letting $\ve \to 0$, we obtain \eqref{F max est}, i.e.
\[ F(y_0,t_0) \le \frac{c(n,k_0)}{\d^2} + c(n)c_1^2,\]
and we can finish the proof as in Case 2.
\end{proof}

\section{A priori estimates of $\nabla^2 g$ along the flow}\label{curvature-section-second-order}

In this section we utilize the arguments in the proof of \cref{est nab g} to
obtain a priori estimates of the second derivatives of $g$ and the Riemann curvature tensor along the Ricci de Turck flow. 
\begin{lem}\label{est nab nab g}
	Under the same assumptions as in \cref{est nab g}, there exists a constant $c(n,k_0) > 0$ depending only on $n$ and $k_0$, such that 
	\begin{align} 
	\begin{split}
	|\wt{\n}^2 g| (x,t) \le	 c(n,k_0) \left( \frac{1}{\d^2} + c_1^2 + \frac{c_2^{1/3}}{\d^{2/3}} + c_2^{1/3} c_1^{2/3}  \right)
	\end{split} 
	\end{align}
	for all $(x,t) \in B(x_0,\g+\d/3) \times [0,T]$, where 
	\[ c_1 = \sup_{x \in B(x_0,\g+3\d/4)} |\wt{\n} \wt{\Rm}|(x), \qquad c_2 = \sup_{x \in B(x_0,\g+3\d/4)} |\wt{\n^2} \wt{\Rm}|(x).  \] 
\end{lem}


\noindent We will prove this result below and first note an immediate consequence: Assuming additionally that $|\wt{\n} \wt{\Rm}| = \mathcal{O}(\rho^{-1})$ and $|\wt{\n}^2 \wt{\Rm}| = \mathcal{O}(\rho^{-2})$, with $\rho > 0$ being the distance to the singularity, we obtain $|\wt{\n}^2 g| = \mathcal{O}(\rho^{-2})$ and 
$|\Rm| = \mathcal{O}(\rho^{-2})$ uniformly in $t \in [0,T]$.

\begin{cor}\label{est nab nab g cor}
	Let $(M,\wt{g})$ be a (possibly incomplete) manifold. Fix $0 < T < \infty$ and let $g(x,t)$ be a smooth solution of 
	the initial value problem
	\begin{equation*}
	\begin{split}
	\frac{\p}{\p t} g_{ij}(x,t) = (-2 \Ric_{ij} + \n_i V_j + \n_j V_i)(x,t)&, \quad  (x,t) \in \, M \times [0,T],  \\
	g(x,0) = \wt{g}(x)&, \quad x \in \, M,
	\end{split}
	\end{equation*}
	where $V^i = g^{jk}(\G^i_{jk} - \wt{\G}^i_{jk})$ is the de Turck vector field. We assume that
	\begin{equation*}
	(1-\ve(n)) \wt{g}(x) \le g(x,t) \le (1+\ve(n)) \wt{g}(x) 
	\end{equation*}
	for $\ve(n) > 0$ sufficiently small, only depending on $n$, and for all $(x,t) \in M \times [0,T]$.	Also assume that 
	\[ |\wt{Rm}|^2 \le k_0 \]
	for some constant $k_0 > 0$, and that there exists a constant $C > 0$ such that 
	for all $x \in M$, $0 < \rho \le 1$
	\[ |\wt{\n}\wt{\Rm}|(x) \le \frac{C}{\rho}, \quad  |\wt{\n}^2\wt{\Rm}|(x) \le \frac{C}{\rho^2} \]
	whenever $B(x,\rho-r)$ is relatively compact for all $r > 0$. Then there exists a constant $C' > 0$ such that for all $x \in M$, $t \in [0,T]$, $0 < \rho \le 1$
	\[ |\wt{\n}^2g|(x,t) \le \frac{C'}{\rho^2}, \ |\Rm|(x,t) \le \frac{C'}{\rho^2} \]
	whenever $B(x,\rho-r)$ is relatively compact for all $r > 0$.
\end{cor}

\begin{proof}[Proof of \cref{est nab nab g cor}]
	Let $x_0 \in M$ and $\rho \le 1$ such that $B(x_0,\rho-r) \subset M$ relatively compact 
	for all $r > 0$. Then by \cref{est nab nab g} (choosing $\g,\d$ equal to $\rho/2$)
	\begin{align} 
	\begin{split}
	|\wt{\n}^2 g| (x,t) \le	 c(n,k_0) \left( \frac{1}{\rho^2} + c_1^2 + \frac{c_2^{1/3}}{\rho^{2/3}} + c_2^{1/3} c_1^{2/3}  \right)
	\end{split} 
	\end{align}
	with the constants estimated by 
	\[ c_1 = \sup_{x \in B(x_0,7\rho/8)} |\wt{\n}\wt{\Rm}|(x) \le \frac{8C}{\rho}, \quad
	c_2 = \sup_{x \in B(x_0,7\rho/8)} |\wt{\n^2} \wt{\Rm}|(x) \le \frac{8\widehat{C}}{\rho^2}, \]
	since for all $x \in B(x_0,7\rho/8)$ we have that $B(x,\rho/8-r) \subset M$ relatively compact for all $r > 0$. 
	The estimate of the Riemannian curvature tensor follows from this, \cref{est nab g cor} and 
	\[ \Rm = \wt{\Rm} * \wt{g}^{-1} * g + \wt{\n}^2 g + g^{-1} * \wt{\n} g * \wt{\n} g,  \]
	see \cite[p. 276, formula (83)]{shi}.
	\end{proof}

\begin{proof}[Proof of \cref{est nab nab g}]
In the following all estimates and inequalities are supposed to hold on $B(x_0,\g+\d/2) \times [0,T]$, when nothing else is mentioned.
Differentiating the equation for the metric $g$ from \cite[Lemma 2.1]{shi} $m$ times we obtain
\begin{align}
\begin{split}
\frac{\p}{\p t} \wt{\n}^m g = & g^{\a\b} \wt{\n}_\a \wt{\n}_\b \wt{\n}^m g \\
& + \sum_{\substack{0 \le k_1,k_2,\dots,k_{m+2} \le m+1 \\ k_1 + k_2 + \dots + k_{m+2} \le m+2}} \wt{\n}^{k_1} g * \wt{\n}^{k_2} g * \dots * \wt{\n}^{k_{m+2}}g * P_{k_1 k_2 \dots k_{m+2}} \\
& + \sum_{\substack{0 \le l_1,l_2,\dots,l_m,s \le m \\ l_1 + l_2 + \dots + l_m + s = m}} \wt{\n}^s \wt{\Rm} * \wt{\n}^{l_1} g * \wt{\n}^l_2 g * \dots * \wt{\n}^{l_m} g * Q_{l_1 l_2 \dots l_m s},
\end{split}
\end{align}
where $P_{k_1 k_2 \dots k_{m+2}}$ and $Q_{l_1 l_2 \dots l_m s}$ are polynomials of $g, g^{-1}$. Hence
\begin{align} \label{d dt nab m est}
\begin{split}
\frac{\p}{\p t} |\wt{\n}^m g|^2 = & g^{\a\b} \wt{\n}_\a \wt{\n}_\b |\wt{\n}^m g|^2 -2g^{\a\b} \wt{\n}_\a \wt{\n}^m g \cdot \wt{\n}_\b \wt{\n}^m g \\
& + \sum_{\substack{0 \le k_1,k_2,\dots,k_{m+2} \le m+1 \\ k_1 + k_2 + \dots + k_{m+2} \le m+2}} \wt{\n}^{k_1} g * \wt{\n}^{k_2} g * \dots * \wt{\n}^{k_{m+2}}g * \wt{\n}^m g * P_{k_1 k_2 \dots k_{m+2}} \\
& + \sum_{\substack{0 \le l_1,l_2,\dots,l_m,s \le m \\ l_1 + l_2 + \dots + l_m + s = m}} \wt{\n}^s \wt{\Rm} * \wt{\n}^{l_1} g * \wt{\n}^l_2 g * \dots * \wt{\n}^{l_m} g * \wt{\n}^m g * Q_{l_1 l_2 \dots l_m s},
\end{split}
\end{align}
For $m=2$ this gives, together with $2 g^{\a\b} \wt{\n}_\a \wt{\n}^2 g \cdot \wt{\n}_\b \wt{\n}^2 g \ge |\wt{\n}^3 g|^2
$,
\begin{align}
\begin{split}
\frac{\p}{\p t} |\wt{\n}^2 g|^2 \le & g^{\a\b} \wt{\n}_\a \wt{\n}_\b |\wt{\n}^2 g|^2 -|\wt{\n}^3 g|^2 \\
& + c(n) |\wt{\n}^2 g| (|\wt{\n}^3 g| |\wt{\n} g| + |\wt{\n}^2 g|^2 + |\wt{\n}^2 g| |\wt{\n}g|^2 + |\wt{\n}g|^4) \\
& + c(n) |\wt{\n}^2 g| (|\wt{\n}^2 \wt{\Rm}| + |\wt{\n}\wt{\Rm}| |\wt{\n}g| + |\wt{\Rm}| |\wt{\n}^2 g| + |\wt{\Rm}| |\wt{\n}g|^2).
\end{split}
\end{align}
It follows that
\begin{align} \label{norm nab squared g inequ}
\begin{split}
\frac{\p}{\p t} |\wt{\n}^2 g|^2 \le & g^{\a\b} \wt{\n}_\a \wt{\n}_\b |\wt{\n}^2 g|^2 -\frac{1}{2} |\wt{\n}^3 g|^2 \\
& + c(n)  (|\wt{\n}^2 g|^2 |\wt{\n}g|^2 + |\wt{\n}^2 g|^3  + |\wt{\n}^2 g| |\wt{\n}g|^4 + c_2 |\wt{\n}^2 g| \\
& \hspace{2cm}+ c_1 |\wt{\n}^2 g| |\wt{\n}g| + \sqrt{k_0} |\wt{\n}^2 g|^2 + \sqrt{k_0} |\wt{\n}^2 g| |\wt{\n}g|^2),
\end{split}
\end{align}
on $B(x_0,\g+3\d/4) \times [0,T]$, where 
\[ c_1 = \sup_{x \in B(x_0,\g+3\d/4)} |\wt{\n} \wt{\Rm}|(x), \qquad c_2 = \sup_{x \in B(x_0,\g+3\d/4)} |\wt{\n^2} \wt{\Rm}|(x),  \]
and where we used $|\wt{\n}^2 g| |\wt{\n}^3 g| |\wt{\n} g| \le \frac{1}{2}|\wt{\n}^3 g|^2 + \frac{1}{2}|\wt{\n}^2 g|^2 |\wt{\n} g|^2$. 
From \eqref{inbetween 1} and \eqref{inbetween 2} we have
\begin{align}\label{norm nab g inequ 1}
\begin{split}
\frac{\p}{\p t} |\wt{\n} g|^2 \le & g^{\a\b} \wt{\n}_\a \wt{\n}_\b |\wt{\n} g|^2 -\frac{1}{2} |\wt{\n}^2 g|^2 + c(n,k_0) |\wt{\n} g|^2 + c(n) c_1 |\wt{\n}g| \\ 
& + 3200n^{10} |\wt{\n} g|^4.
\end{split}
\end{align}
Now as in \cite[Proof of Lemma 4.2, p.256 (80)]{shi} let
\begin{align}\label{phi2-definition}
\psi(x,t) = (a+|\wt{\n}g|^2) |\wt{\n}^2 g|^2,
\end{align}
where $a > 0$ is a constant which is chosen later. Then
\begin{align}
\begin{split}
\left(\frac{\p}{\p t} - g^{\a\b} \wt{\n}_\a \wt{\n}_\b \right)\psi = & \left(\frac{\p}{\p t} - g^{\a\b} \wt{\n}_\a \wt{\n}_\b \right) ( a+|\wt{\n}g|^2) \cdot |\wt{\n}^2 g|^2 \\
& + (a + |\wt{\n}g|^2) \left(\frac{\p}{\p t} - g^{\a\b} \wt{\n}_\a \wt{\n}_\b \right) |\wt{\n}^2 g|^2 \\
& - 2 g^{\a\b} \wt{\n}_\a |\wt{\n}g|^2 \wt{\n}_\b |\wt{\n}^2 g|^2.
\end{split}
\end{align}
We proceed as before in Lemmas \ref{Shi Lemma 4.1} and \ref{est nab g} along the following steps.
\medskip

\begin{enumerate}
\item Step 1: Derive an evolution inequality for $\psi$. 
\smallskip

\item Step 2: Estimate $\wt{\n} \wt{\n} \xi$ from below.
\smallskip

\item Step 3: Estimate $\xi \psi$ from above and conclude the proof.
\medskip

\end{enumerate}

\noindent \textbf{Step 1: Derive an evolution inequality for $\psi$.}
\medskip

\noindent Together with \eqref{norm nab squared g inequ} and \eqref{norm nab g inequ 1} we obtain
\begin{align}
\begin{split}
\frac{\p}{\p t} \psi \le & g^{\a\b} \wt{\n}_\a \wt{\n}_\b \psi -\frac{1}{2} |\wt{\n}^2 g|^4 + c(n,k_0) |\wt{\n}g|^2 |\wt{\n}^2 g|^2 + c(n)c_1 |\wt{\n}g| |\wt{\n}^2 g|^2 \\
& + c(n) |\wt{\n}g|^4 |\wt{\n}^2 g|^2 -\frac{1}{2} (a + |\wt{\n}g|^2) |\wt{\n}^3 g|^2 \\
& + c(n,k_0) (a + |\wt{\n}g|^2) (|\wt{\n}^2 g|^2 |\wt{\n}g|^2 + |\wt{\n}^2 g|^3  + |\wt{\n}^2 g| |\wt{\n}g|^4 + c_2 |\wt{\n}^2 g| \\
& + c_1 |\wt{\n}^2 g| |\wt{\n}g| + |\wt{\n}^2 g|^2 + |\wt{\n}^2 g| |\wt{\n}g|^2) \\
& - 2 g^{\a\b} \wt{\n}_\a |\wt{\n}g|^2 \wt{\n}_\b |\wt{\n}^2 g|^2.
\end{split}
\end{align}
We estimate the last term as
\begin{align}
\begin{split}
- 2 g^{\a\b} \wt{\n}_\a |\wt{\n}g|^2 \wt{\n}_\b |\wt{\n}^2 g|^2 & \le 16 |\wt{\n}g| |\wt{\n}^2 g|^2 |\wt{\n}^3 g| \\
& \le 16 C_1 |\wt{\n}^2 g|^2 |\wt{\n}^3 g| \\
& \le \frac{1}{2}a |\wt{\n}^3 g|^2 + \frac{1}{2a} \cdot 256 C_1^2 |\wt{\n}^2 g|^4 \\
& = \frac{1}{2}a |\wt{\n}^3 g|^2 + \frac{1}{4} |\wt{\n}^2 g|^4,
\end{split}
\end{align}
where $C_1 := \frac{c(n,k_0)}{\d} + c(n)c_1$ is the bound on $|\wt{\n}g|$ from \eqref{nab g est lemma} and we chose $a = 512 C_1^2$. This gives
\begin{align}
\begin{split}
\frac{\p}{\p t} \psi \le & g^{\a\b} \wt{\n}_\a \wt{\n}_\b \psi -\frac{1}{4} |\wt{\n}^2 g|^4 + c(n,k_0) C_1^2 |\wt{\n}^2 g|^2 + c(n) c_1 C_1 |\wt{\n}^2 g|^2 + c(n)C_1^4 |\wt{\n}^2 g|^2 \\
& + c(n,k_0) C_1^2 (C_1^2 |\wt{\n}^2 g|^2 + |\wt{\n}^2 g|^3 + C_1^4 |\wt{\n}^2 g| + c_2 |\wt{\n}^2 g| + c_1 C_1 |\wt{\n}^2 g| \\
& \hspace{3cm}+ |\wt{\n}^2 g|^2  + C_1^2 |\wt{\n}^2 g|). 
\end{split}
\end{align}
Now by definition of $\psi$ we have
\begin{align}
|\wt{\n}^2 g|^2 = \frac{\psi}{a+|\wt{\n}g|^2} \le \frac{\psi}{a} = \frac{\psi}{512 C_1^2}
\end{align}
and 
\begin{align}
|\wt{\n}^2 g|^2 = \frac{\psi}{a+|\wt{\n}g|^2} \ge \frac{\psi}{a+C_1^2} = \frac{\psi}{513C_1^2}
\end{align}
This yields
\begin{align}\label{psi inequ 1}
\begin{split}
\frac{\p}{\p t}\psi & \le g^{\a\b} \wt{\n}_\a \wt{\n}_\b \psi -\frac{1}{4} \frac{\psi^2}{513^2 C_1^4} + c(n,k_0) (\psi + \frac{c_1}{C_1} \psi + C_1^2 \psi) \\
&+ c(n,k_0) (C_1^2 \psi + \frac{\psi^{3/2}}{C_1} + C_1^5 \psi^{1/2} + c_2 C_1 \psi^{1/2} + c_1 C_1^2 \psi^{1/2} + \psi + C_1^3 \psi^{1/2}) \\
& = g^{\a\b} \wt{\n}_\a \wt{\n}_\b \psi -\frac{1}{4} \frac{\psi^2}{513^2 C_1^4} \\
&+ c(n,k_0) \left(\frac{\psi^{3/2}}{C_1} + (C_1^2 + \frac{c_1}{C_1} + 1) \psi + (C_1^5 + C_1^3 + c_1 C_1^2 + c_2 C_1) \psi^{1/2} \right) \\
& \le g^{\a\b} \wt{\n}_\a \wt{\n}_\b \psi -\frac{1}{4} \frac{\psi^2}{513^2 C_1^4} + c(n,k_0) \left(\frac{\psi^{3/2}}{C_1} + C_1^2 \psi + (C_1^5 + c_2 C_1) \psi^{1/2} \right), 
\end{split}
\end{align}
where in the last step we used that $C_1 = \frac{c(n,k_0)}{\d} + c(n)c_1 \ge c_1$ and that, 
assuming $c(n,k_0) \ge 1$, $C_1 \ge 1$. \medskip

\noindent \textbf{Step 2: Estimate $\wt{\n} \wt{\n} \xi$ from below.}
\medskip

\noindent Now let $\eta \in C^\infty(\r)$ be the cutoff function as before and define the cutoff function $\xi \in C(M)$ as 
\begin{align}
\begin{split}
\xi(x) = \eta \left(\frac{d_{\wt{g}}(x,x_0) - (\g+\d/3)}{\d/12}\right),
\end{split}
\end{align}
where $d_{\wt{g}}$ denotes the distance function with respect to the metric $\wt{g}$. Then we have 
\begin{align}
\begin{split}
& \xi(x) = 1, \qquad x \in B(x_0,\g+\d/3), \\
& \xi(x) = 0, \qquad x \in M \backslash B(x_0,\g+5\d/12), \\
& 0 \le \xi(x) \le 1, \qquad x \in M.
\end{split}
\end{align} 
If $d_{\wt{g}}(\cdot,x_0)$ is smooth in a neighborhood of a point $x$, we have by a calculation analogous to \eqref{inbetween 3}
\begin{align}
|\wt{\n} \xi(x)|^2 \le \frac{2304}{\d^2} \xi(x)
\end{align}
and 
\begin{align}\label{est nab squared xi}
\begin{split}
\wt{\n}_\a \wt{\n}_\b \xi(x) \ge -\left( \frac{1152}{\d^2} + \frac{48}{\d} \sqrt[4]{k_0} \coth \left(\sqrt[4]{k_0} d_{\wt{g}}(x,x_0) \right) \right) \wt{g}_{\a\b}(x).
\end{split}
\end{align} \ \medskip

\noindent \textbf{Step 3: Estimate $\xi \psi$ from above and conclude the proof.}
\medskip

Let
\[ F(x,t) := \xi(x)\psi(x,t), \qquad (x,t) \in B(x_0,\g+\d) \times [0,T]. \]
Since $|\wt{\n}^2 g|^2(x,0) = 0$, we have
\begin{align}
F(x,0) = 0, \qquad x \in B(x_0,\g+\d).
\end{align}
Since $\xi(x) = 0$ for $x \in B(x_0,\g+\d) \backslash B(x_0,\g+\frac{5}{12}\d)$, it follows that
\begin{align}
F(x,t) = 0, \qquad (x,t) \in B(x_0,\g+\d) \backslash B(x_0,\g+\frac{5}{12}\d) \times [0,T].
\end{align}
Thus there exists a point $(y_0,t_0) \in B(x_0,\g+\frac{5}{12}\d) \times [0,T]$ with $t_0 > 0$ such that
\begin{align}
F(y_0,t_0) =  \max \ \left\{F(x,t) \mid (x,t) \in B(x_0,\g+\d) \times [0,T] \right\}
\end{align}
unless $F \equiv 0$ on $B(x_0,\g+\d) \times [0,T]$. \medskip

Next, as previously in Lemma \ref{est nab g}, 
we distinguish three cases, first case where $\xi \equiv 1$ in a neighborhood of $y_0$, 
second case where $\xi$ is not identically $1$, but smooth in a neighborhood of $y_0$, 
and third case, where $\xi$ is not smooth and a trick needs to be applied. \medskip

\noindent \underline{\textbf{Case 1. $y_0 \in B(x_0,\g+\frac{5}{12}\d)$}} \medskip

\noindent Then $\xi \equiv 1$ in a neighborhood of $y_0$, such that $F = \psi$ near $(y_0,t_0)$, and we have 
\begin{align}
0 & \le (\frac{\p}{\p t} - g^{\a\b}\wt{\n}_\a \wt{\n}_\b) \psi(y_0,t_0) \\
& \le -\frac{1}{4} \frac{\psi^2(y_0,t_0)}{513^2 C_1^4} + c(n,k_0) \left(\frac{\psi^{3/2}}{C_1} + C_1^2 \psi + (C_1^5 + c_2 C_1) \psi^{1/2} \right)(y_0,t_0)
\end{align} 
and thus 
\begin{align}
\frac{1}{4}\frac{F^2(y_0,t_0)}{513^2 C_1^4} \le c(n,k_0) \left(\frac{F^{3/2}}{C_1} + C_1^2 F + (C_1^5 + c_2 C_1) F^{1/2} \right)(y_0,t_0)
\end{align}
which is a better estimate than the one below in Case 2, and thus Case 1 follows from Case 2. \\

\noindent \underline{\textbf{Case 2. $y_0 \notin B(x_0,\g+\frac{5}{12}\d)$ and $y_0$ is not in the cut locus of $x_0$}}\medskip

\noindent Then the distance function $d_{\wt{g}}(\cdot,x_0)$ is smooth in a neighborhood of $y_0$ and it follows that
\begin{align}\label{F inequ at max}
\begin{split}
& 0 \le \frac{\p F}{\p t}(y_0,t_0) = \xi(y_0) \frac{\p \psi}{\p t}(y_0,t_0), \\
& 0 = \wt{\n}_\a F(y_0,t_0) = (\xi \wt{\n}_\a \psi + \psi \wt{\n}_\a \xi)(y_0,t_0), \\
& 0 \ge g^{\a\b} \wt{\n}_\a \wt{\n}_\b F(y_0,t_0) = (\xi g^{\a\b} \wt{\n}_\a \wt{\n}_\b \psi + \psi g^{\a\b} \wt{\n}_\a \wt{\n}_\b \xi 
\\ &+ 2 g^{\a\b} \wt{\n}_\a \xi \wt{\n}_\b \psi)(y_0,t_0).
\end{split}
\end{align}
Together with \eqref{psi inequ 1} we obtain at the point $(y_0,t_0)$
\begin{align}
\begin{split}
\frac{1}{4}\frac{\psi^2}{513^2 C_1^4}\xi & \le \xi g^{\a\b}\wt{\n}_\a \wt{\n}_\b \psi + \xi c(n,k_0) \left(\frac{\psi^{3/2}}{C_1} + C_1^2 \psi + (C_1^5 + c_2 C_1) \psi^{1/2} \right)\\
& \le -\psi g^{\a\b} \wt{\n}_\a \wt{\n}_\b \xi - 2 g^{\a\b} \wt{\n}_\a \xi \wt{\n}_\b \psi \\
& \hspace{0.5cm}+ \xi c(n,k_0) \left(\frac{\psi^{3/2}}{C_1} + C_1^2 \psi + (C_1^5 + c_2 C_1) \psi^{1/2} \right). 
\end{split}
\end{align}
From \eqref{F inequ at max} we have at $(y_0,t_0)$
\begin{align}
- 2 g^{\a\b} \wt{\n}_\a \xi \wt{\n}_\b \psi  = \frac{2\psi}{\xi} g^{\a\b} \wt{\n}_\a \xi \wt{\n}_\b \xi \le \frac{9216}{\d^2} \psi,
\end{align}
where the last inequality follows from \eqref{equiv met 2}. Furthermore, from \eqref{est nab squared xi} and an estimate analogous to \eqref{est coth} we obtain at $(y_0,t_0)$
\begin{align}
-\psi g^{\a\b} \wt{\n}_\a \wt{\n}_\b \xi \le \frac{c(n,k_0)}{\d^2} \psi.
\end{align}
This yields the following intermediate inequality
\begin{align}
\begin{split}
\frac{1}{4}\frac{\psi^2}{513^2 C_1^4}\xi \le & \frac{c(n,k_0)}{\d^2} \psi + \xi c(n,k_0) \left(\frac{\psi^{3/2}}{C_1} + C_1^2 \psi + (C_1^5 + c_2 C_1) \psi^{1/2} \right). 
\end{split}
\end{align}
Multiplying this inequality with $\xi$, using $0 \le \xi \le 1$ and adjusting the constants $c(n,k_0)>0$ 
appropriately, we obtain
\begin{align}\label{second-derivatie-estimate1}
\begin{split}
\frac{F(y_0,t_0)^2}{C_1^4} & \le c(n,k_0) \left(\frac{F(y_0,t_0)^{3/2}}{C_1} + C_1^2 F(y_0,t_0) + (C_1^5 + c_2 C_1) F(y_0,t_0)^{1/2} \right)
\\ &+ \frac{c(n,k_0)}{\d^2} F(y_0,t_0).
\end{split}
\end{align}
Now we use the following elementary estimate: If $x \ge 0$ satisfies
\begin{equation} \label{elem est} 
x^2 \le a x^{3/2} + bx + cx^{1/2}
\end{equation}
with constants $a,b,c \ge 0$, then 
\[ x \le \max\{a^2,b,c^{2/3}\}.\] 
This reduces \eqref{second-derivatie-estimate1} to the following estimate
\begin{align}
\begin{split}
F(y_0,t_0) & \le c(n,k_0) \left(C_1^6 + \frac{C_1^4}{\d^2} + (C_1^9 + c_2 C_1^5)^{2/3} \right)\\
& \le c(n,k_0) \left(C_1^6 + (C_1^9 + c_2 C_1^5)^{2/3} \right),
\end{split}
\end{align}
since $\frac{1}{\d} \le C_1$. It follows that for all $(x,t) \in B(x_0,\g+\d) \times [0,T]$
\[ F(x,t) \le F(y_0,t_0) \le c(n,k_0) \left(C_1^6 + (C_1^9 + c_2 C_1^5)^{2/3} \right). \]
Since $F(x,t) = \xi(x) \psi(x,t)$ and $\xi(x) = 1$ for $x \in B(x_0,\g+\d/3)$, we obtain
\[ \psi(x,t) \le c(n,k_0) \left(C_1^6 + (C_1^9 + c_2 C_1^5)^{2/3} \right) \qquad \forall (x,t) \in B(x_0,\g+\d/3) \times [0,T]. \]
As $\psi(x,t) = (a + |\wt{\n}g|^2) |\wt{\n}^2 g|^2$ and $a = 512 C_1^2$ we have
\[ |\wt{\n}^2 g|^2 (x,t) = \frac{\psi(x,t)}{(a + |\wt{\n}g|^2)(x,t)} \le \frac{\psi(x,t)}{a} \le c(n,k_0) \left(C_1^4 + (C_1^6 + c_2 C_1^2)^{2/3} \right) \]
for all $(x,t) \in B(x_0,\g+\d/3) \times [0,T]$. Thus 
\begin{align} 
\begin{split}
|\wt{\n}^2 g| (x,t) & \le c(n,k_0) \sqrt{C_1^4 + (C_1^6 + c_2 C_1^2)^{2/3}} \\
& \le c(n,k_0) (C_1^2 + (C_1^6 + c_2 C_1^2)^{1/3}) \\
& \le c(n,k_0) (C_1^2 + c_2^{1/3} C_1^{2/3}) \\
& \le c(n,k_0) \left( \left(\frac{1}{\d}+c_1 \right)^2 + c_2^{1/3} \left(\frac{1}{\d}+c_1 \right)^{2/3} \right) \\
& \le c(n,k_0) \left( \frac{1}{\d^2} + c_1^2 + \frac{c_2^{1/3}}{\d^{2/3}} + c_2^{1/3} c_1^{2/3}  \right)
\end{split} 
\end{align}
for all $(x,t) \in B(x_0,\g+\d/3) \times [0,T]$, where we used 
\[ \sqrt[3]{a+b} \le \sqrt[3]{a} + \sqrt[3]{b}, \quad (a+b)^2 \le 2a^2 + 2b^2\] 
for real numbers $a,b \ge 0$. \medskip

\noindent \underline{\textbf{Case 3. $y_0 \notin B(x_0,\g+\frac{5}{12}\d)$ and $y_0$ is in the cut locus of $x_0$}} \medskip

\noindent Then we again apply Calabi's trick, see Case 3 in the Proof of \cref{est nab g}.
\end{proof}

\section{A priori estimates of $\nabla^m g$ along the flow}\label{curvature-section-higher-order}

In this section prove we prove a priori estimates for all higher derivatives of $g$ and the Riemann curvature tensor along the Ricci de Turck flow. We treated the case of the second derivatives $\wt{\n}^2 g$ separately, since the evolution inequality \eqref{norm nab g inequ 1} for $|\wt{\n} g|^2$ which goes into the estimate of the time-derivative of $\psi$ (see the proof of \cref{est nab nab g} above) differs from the 
corresponding one \eqref{est nab^m-1 g proof} below that will be obtained for the higher derivatives $|\wt{\n}^{m-1} g|^2$.  

\begin{lem}\label{est nab^m g}
	Under the same assumptions as in \cref{est nab g} we set for $k, s \in \N_0$
	\begin{equation}
	\begin{split}
	C_k := \sup_{x \in B(x_0,\g+\d/(k+1))} |\wt{\n}^k g|, \quad
	c_s := \sup_{x \in B(x_0,\g+3\d/4)} |\wt{\n}^s \wt{\Rm}|,
	\end{split}
	\end{equation}
	and define for any integer $p \geq 1$ the following constants
	\begin{equation}
	\begin{split}
&\mathcal{K}_p := \sum_{\substack{0 \le k_1,\dots,k_{p+2} \le p-1 \\ k_1 + \dots + k_{p+2} \le p+2}} C_{k_1} \cdots C_{k_{p+2}}, \\
&\mathcal{L}_p := \sum_{\substack{0 \le l_1,\dots,l_p,s \le p-1 \\ l_1 + \dots + l_p + s = p}} c_s C_{l_1} \cdots C_{l_p}. 
	\end{split}
	\end{equation}
	Then we find for $m\geq 3$ and for all $(x,t) \in B(x_0,\g+\frac{\d}{m+1}) \times [0,T]$
	\begin{align} 
	\begin{split}
	 |\wt{\n}^m g|^2 (x,t) \le \max \, \{ A, B\}
	\end{split}
	\end{align}
	where for some constants $c(n,m,k_0), c(n,m) > 0$
	\begin{equation}
	\begin{split}
	&A:= c(n,m,k_0) C_{m-1}^2 \left( \frac{1}{\d^2} + C_1^2 + C_2 + 
	\sqrt{k_0} + \frac{\mathcal{K}_{m-1} + \mathcal{L}_{m-1} + c_{m-1}}{C_{m-1}}    \right), \\
	&B:= c(n,m) \frac{1}{C_{m-1}^2}(C_{m-1}^5 (\mathcal{K}_{m} + \mathcal{L}_{m} + c_{m}))^{2/3}. 
	\end{split}
	\end{equation}
	\end{lem}
We first prove a corollary of that result and later provide the proof of the lemma above. We point out that
with more effort it would be possible to obtain an even more explicit bound of $|\wt{\n}^m g|^2$ analogous to the one in \cref{est nab nab g}, but since our main interest is in the behaviour of the derivatives of the metric and the Riemann curvature tensor when approaching the singular strata, the bound above is sufficient for our purposes.

\begin{cor}\label{cor nabla^m}
	Let $(M,\wt{g})$ be a (possibly incomplete) manifold. Fix $0 < T < \infty$ and let $g(x,t)$ be a smooth solution of 
	the initial value problem
	\begin{equation*}
	\begin{split}
	\frac{\p}{\p t} g_{ij}(x,t) = (-2 \Ric_{ij} + \n_i V_j + \n_j V_i)(x,t)&, \quad  (x,t) \in \, M \times [0,T],  \\
	g(x,0) = \wt{g}(x)&, \quad x \in \, M,
	\end{split}
	\end{equation*}
	where $V^i = g^{jk}(\G^i_{jk} - \wt{\G}^i_{jk})$ is the de Turck vector field. We assume that
	\begin{equation*}
	(1-\ve(n)) \wt{g}(x) \le g(x,t) \le (1+\ve(n)) \wt{g}(x) 
	\end{equation*}
	for $\ve(n) > 0$ sufficiently small, only depending on $n$, and for all $(x,t) \in M \times [0,T]$.	Also assume that 
	\[ |\wt{Rm}|^2 \le k_0 \]
	for some constant $k_0 > 0$, and that for all $m \ge 1$ there exists a constant $C_m > 0$, such that for all $x \in M$, $0 < \rho \le 1$
	\[ |\wt{\n}^m \wt{\Rm}|(x) \le \frac{C}{\rho^m}  \] 
	whenever $B(x,\rho-r)$ is relatively compact for all $r > 0$. Then there exists a constant $C'_m > 0$ such that for all $x \in M$, $t \in [0,T]$, $0 < \rho \le 1$
	\[ |\wt{\n}^m g|(x,t) \le \frac{C'_m}{\rho^m}, \quad |\n^m \Rm|(x,t) \le \frac{C_m'}{\rho^{m+2}} \]
	whenever $B(x,\rho-r)$ is relatively compact for all $r > 0$.
\end{cor}

\begin{proof}
We start with the estimates of the derivatives of the metric $g$. The cases $m=1,2$ have already been proven, so assume that $m \ge 3$. By induction, we can assume that there exists a constant $C' > 0$ such that for all $k = 1,\dots,m-1$, $(x,t) \in M \times [0,T]$, $\rho \le 1$, $r > 0$ 
\[ |\wt{\n}^k g|(x,t) \le \frac{C'}{\rho^k} \] 
whenever $B(x,\rho-r)\subset M$ is relatively compact. 
Let $x_0 \in M$ and $\rho \le 1$ such that $B(x_0,\rho-r) \subset M$ relatively compact 
for all $r > 0$. Then by \cref{est nab nab g} (choosing $\g,\d$ equal to $\rho/2$)
\begin{align} 
	\begin{split}
	 |\wt{\n}^m g|^2 (x,t) \le \max \, \{ A, B\}
	\end{split}
	\end{align}
	for all $t \in [0,T]$. Recall the explicit form of $A$ and $B$
	\begin{equation}
	\begin{split}
	&A= c(n,m,k_0) C_{m-1}^2 \left( \frac{1}{\d^2} + C_1^2 + C_2 + 
	\sqrt{k_0} + \frac{\mathcal{K}_{m-1} + \mathcal{L}_{m-1} + c_{m-1}}{C_{m-1}}    \right), \\
	&B= c(n,m) \frac{1}{C_{m-1}^2}(C_{m-1}^5 (\mathcal{K}_{m} + \mathcal{L}_{m} + c_{m}))^{2/3}. 
	\end{split}
	\end{equation}
The individual constants can be estimated as follows:
\begin{align}
\begin{split}
C_k = \sup_{x \in B(x_0,\rho/2+\rho/2/(k+1))} |\wt{\n}^k g| \le \sup_{x \in B(x_0,3\rho/4)} |\wt{\n}^k g| \le \frac{4^k C'}{\rho^k}
\end{split}
\end{align}
for $k = 1,\dots,m-1$, since for all $x \in B(x_0,3\rho/4)$ we have that $B(x,\rho/4-r) \subset M$ relatively compact for all $r > 0$,
\begin{align}
\begin{split}
c_s = \sup_{x \in B(x_0,\rho/2+3\rho/2/4)} |\wt{\n}^s \wt{\Rm}| = \sup_{x \in B(x_0,7\rho/8)} |\wt{\n}^s \wt{\Rm}| \le \frac{8^s C}{\rho^s}
\end{split}
\end{align}
for $s = 1,\dots,m$, since for all $x \in B(x_0,7\rho/8)$ we have that $B(x,\rho/8-r) \subset M$ relatively compact for all $r > 0$. Thus
\begin{align}
\begin{split}
\mathcal{K}_m \le \frac{C}{\rho^{m+2}}, \quad \mathcal{L}_m \le \frac{C}{\rho^m}, \quad \mathcal{K}_{m-1} \le \frac{C}{\rho^{m+1}}, \quad \mathcal{L}_{m-1} \le \frac{C}{\rho^{m-1}}
\end{split}
\end{align}
with the constant $C > 0$ only depending on $m$ and the constants $C', C$ from above. 
Plugging this in gives 
\[ |\wt{\n}^m g|^2 (x_0,t) \le \frac{C}{\rho^{2m}}\]
for all $t \in [0,T]$, with $C > 0$ only depending on $m, n, k_0$ and the constants $C', C$ from above.	This completes the proof for the derivatives of the metric. \medskip

To estimate the derivatives of the curvature tensor, we start by the following general identities for any (say $(1,2)$-tensor) $A$
\[ \n_l A^i_{jk} = \frac{\p}{\p x^l} A^i_{jk} + A^m_{jk} \G^i_{ml} - A^i_{mk} \G^m_{jl} - A^i_{jm} \G^m_{kl},\]	
\[ \wt{\n}_l A^i_{jk} = \frac{\p}{\p x^l} A^i_{jk} + A^m_{jk} \wt{\G}^i_{ml} - A^i_{mk} \wt{\G}^m_{jl} - A^i_{jm} \wt{\G}^m_{kl}.\]
Thus $\n$ and $\wt{\n}$, acting on $(1,2)$-tensors, differ by the following expression
\begin{equation}\label{A identity}
\n_l A^i_{jk} = \wt{\n}_l A^i_{jk} + A^m_{jk} (\G^i_{ml} - \wt{\G}^i_{ml}) - A^i_{mk} (\G^m_{jl} - \wt{\G}^m_{jl}) - A^i_{jm} (\G^m_{kl} - \wt{\G}^m_{kl}). 
\end{equation}
In normal coordinates at a point $p \in M$ with respect to the metric $\wt{g}$ we have $\wt{\G}^k_{ij} = 0$ and $\frac{\p}{\p x^i} g_{jk} = \wt{\n}_i g_{jk}$ at the point $p$, such that
\[ \G^k_{ij} - \wt{\G}^k_{ij} = \frac{1}{2} g^{km} (\wt{\n}_j g_{im} + \wt{\n}_i g_{jm} - \wt{\n}_m g_{ij}) \]
at $p$. But since this is an identity of tensors, it actually holds for all points in any coordinate system. Using the $*$-notation we can write this shorter as 
\[ \G - \wt{\G} = g^{-1} * \wt{\n} g.\]
Hence \eqref{A identity} takes the form
\[ \n A = \wt{\n} A + A * g^{-1} * \wt{\n} g. \]	
By induction, together with the product rule 
\[ \wt{\n}(A*B) = \wt{\n}A * B + A * \wt{\n}B, \] and the covariant derivative of the inverse metric tensor given by 
\[ \wt{\n}(g^{-1}) = g^{-1} * g^{-1} * \wt{\n}g, \]
we obtain for all $k \ge 1$
\begin{equation}\label{nabla^k A est}
\n^k A = \sum_{\substack{0 \le k_1,\dots,k_r \le k \\ k_1 + \dots + k_r = k}} \wt{\n}^{k_1} A * \wt{\n}^{k_2} g * \dots * \wt{\n}^{k_r} g * P_{k_1 \dots k_r}, 
\end{equation}
where $P_{k_1 \dots k_r}$ is a polynomial in $g^{-1}$.
Now from the identity for the Riemann curvature tensor 
\[ \Rm = \wt{\Rm} * \wt{g}^{-1} * g + \wt{\n}^2 g + g^{-1} * \wt{\n} g * \wt{\n} g, \]
see \cite[p. 276, formula (83)]{shi}, we obtain by induction for all $k \ge 1$
\begin{align} \label{wt nabla^k Rm est} 
\begin{split}
\wt{\n}^k \Rm = & \sum_{\substack{0 \le s,k_1,\dots,k_r \le k \\ s+k_1 + \dots + k_r = k}} \wt{\n}^s \wt{\Rm} * \wt{\n}^{k_1} g * \dots * \wt{\n}^{k_r} g * Q_{sk_1\dots k_r} \\
& + \sum_{\substack{0 \le l_1,\dots,l_s \le k+2 \\ l_1 + \dots + l_s = k+2}} \wt{\n}^{l_1} g * \dots * \wt{\n}^{l_s} g * R_{l_1\dots l_s},
\end{split}  
\end{align}
where $Q,R$ are polynomials in $g,g^{-1}$ and $\wt{g}^{-1}$. 
Plugging \eqref{wt nabla^k Rm est} into \eqref{nabla^k A est} gives
\begin{align} 
\begin{split}
\n^k \Rm = & \sum_{\substack{0 \le s,k_1,\dots,k_r \le k \\ s+k_1 + \dots + k_r = k}} \wt{\n}^s \wt{\Rm} * \wt{\n}^{k_1} g * \dots * \wt{\n}^{k_r} g * S_{sk_1\dots k_r} \\
& + \sum_{\substack{0 \le l_1,\dots,l_s \le k+2 \\ l_1 + \dots + l_s = k+2}} \wt{\n}^{l_1} g * \dots * \wt{\n}^{l_s} g * T_{l_1\dots l_s},  
\end{split}
\end{align}
where $S,T$ are polynomials in $g,g^{-1}$ and $\wt{g}^{-1}$, and thus 
\begin{align*}
\begin{split} 
& |\n^k \Rm| \le \\
& C(n,k) \left( \sum_{\substack{0 \le s,k_1,\dots,k_r \le k \\ s+k_1 + \dots + k_r = k}} |\wt{\n}^s \wt{\Rm}| \cdot |\wt{\n}^{k_1} g| \cdots  |\wt{\n}^{k_r} g| + \sum_{\substack{0 \le l_1,\dots,l_s \le k+2 \\ l_1 + \dots + l_s = k+2}} |\wt{\n}^{l_1} g| \cdots |\wt{\n}^{l_s} g| \right).  
\end{split}
\end{align*}
Now the claim follows from the estimates of the derivatives of $g$.
\end{proof}

\begin{proof}[Proof of \cref{est nab^m g}]
Let $m \ge 2$. From \eqref{d dt nab m est} and since 
$$
2 g^{\a\b} \wt{\n}_\a \wt{\n}^m g \cdot \wt{\n}_\b \wt{\n}^m g \ge |\wt{\n}^{m+1} g|^2,
$$ we have the following differential inequality
\begin{align}
\begin{split}
\frac{\p}{\p t} |\wt{\n}^m g|^2 & \le g^{\a\b} \wt{\n}_\a \wt{\n}_\b |\wt{\n}^m g|^2 - |\wt{\n}^{m+1} g|^2 \\
& \hspace{0.5cm} + c(n) |\wt{\n}^m g| \sum_{\substack{0 \le k_1,\dots,k_{m+2} \le m+1 \\ k_1 + \dots + k_{m+2} \le m+2}} |\wt{\n}^{k_1} g| \cdots |\wt{\n}^{k_{m+2}} g| \\
& \hspace{0.5cm} + c(n) |\wt{\n}^m g| \sum_{\substack{0 \le l_1,\dots,l_m,s \le m \\ l_1 + \dots + l_m + s = m}} |\wt{\n}^s \wt{\Rm}| |\wt{\n}^{l_1} g| \cdots |\wt{\n}^{l_m} g| \\
& \le g^{\a\b} \wt{\n}_\a \wt{\n}_\b |\wt{\n}^m g|^2 - |\wt{\n}^{m+1}|^2 \\ 
& \hspace{0.5cm} + c(n,m) |\wt{\n}^m g| \cdot  [ |\wt{\n}^{m+1} g| |\wt{\n} g| + |\wt{\n}^m g| |\wt{\n}^2 g| + |\wt{\n}^m g| |\wt{\n} g|^2 \\
& \hspace{0.5cm} + \sum_{\substack{0 \le k_1,\dots,k_{m+2} \le m-1 \\ k_1 + \dots + k_{m+2} \le m+2}} C_{k_1} \cdots C_{k_{m+2}}  ] \\
&  \hspace{0.5cm} + c(n,m) |\wt{\n}^m g| \cdot [ |\wt{\Rm}| |\wt{\n}^m g| + |\wt{\n}^m \wt{\Rm}| \\
& \hspace{0.5cm} + \sum_{\substack{0 \le l_1,\dots,l_m,s \le m-1 \\ l_1 + \dots + l_m + s = m}} c_s C_{l_1} \cdots C_{l_m}   ]
\end{split}
\end{align}	
on $B(x_0,\g+\d/m) \times [0,T]$, where we have set as before
\[ c_s := \sup_{x \in B(x_0,\g+3\d/4)} |\wt{\n}^s \wt{\Rm}|, \qquad C_k := \sup_{x \in B(x_0,\g+\d/(k+1))} |\wt{\n}^k g|. \] The following estimates also hold on $B(x_0,\g+\d/m) \times [0,T]$, when nothing else is mentioned. With the abbreviations 
\[ \mathcal{K}_m := \sum_{\substack{0 \le k_1,\dots,k_{m+2} \le m-1 \\ k_1 + \dots + k_{m+2} \le m+2}} C_{k_1} \cdots C_{k_{m+2}}, \qquad \mathcal{L}_m := \sum_{\substack{0 \le l_1,\dots,l_m,s \le m-1 \\ l_1 + \dots + l_m + s = m}} c_s C_{l_1} \cdots C_{l_m} \] and using $|\wt{\n}^m g| |\wt{\n}^{m+1} g| |\wt{\n} g| \le \frac{1}{2}|\wt{\n}^{m+1} g|^2 + \frac{1}{2} |\wt{\n}^m g|^2 |\wt{\n} g|^2$ we obtain
\begin{align} \label{est nab^m g proof}
\begin{split}
\frac{\p}{\p t} |\wt{\n}^m g|^2 & \le g^{\a\b} \wt{\n}_\a \wt{\n}_\b |\wt{\n}^m g|^2 - \frac{1}{2} |\wt{\n}^{m+1}|^2 \\
& \hspace{0.5cm} + c(n,m) \cdot [ |\wt{\n} g|^2 |\wt{\n}^m g|^2 + |\wt{\n}^2 g| |\wt{\n}^m g|^2 + \mathcal{K}_m \cdot |\wt{\n}^m g| \\
& \hspace{0.5cm} + \sqrt{k_0} |\wt{\n}^m g|^2 + c_m |\wt{\n}^m g| + \mathcal{L}_m \cdot |\wt{\n}^m g|  ] \\
& \le g^{\a\b} \wt{\n}_\a \wt{\n}_\b |\wt{\n}^m g|^2 - \frac{1}{2} |\wt{\n}^{m+1}|^2 \\
& \hspace{0.5cm} + c(n,m) \cdot \left( |\wt{\n}^m g|^2 (C_1^2 + C_2 + \sqrt{k_0}) \right) \\
& \hspace{0.5cm} + c(n,m) \cdot  \left(|\wt{\n}^m g| (\mathcal{K}_m + \mathcal{L}_m + c_m) \right).  
\end{split}
\end{align}
Assume from now on that $m \ge 3$. Then we can replace $m$ by $m-1$ and obtain 
\begin{align} \label{est nab^m-1 g proof}
\begin{split}
\frac{\p}{\p t} |\wt{\n}^{m-1} g|^2 & \le g^{\a\b} \wt{\n}_\a \wt{\n}_\b |\wt{\n}^m g|^2 - \frac{1}{2} |\wt{\n}^{m}|^2 \\
& + c(n,m-1) \cdot \left( C_{m-1}^2 (C_1^2 + C_2 + \sqrt{k_0}) \right. \\
& + \left. C_{m-1} (\mathcal{K}_{m-1} + \mathcal{L}_{m-1} + c_{m-1}) \right).  
\end{split}
\end{align}
We define similar to \eqref{phi2-definition}
\[ \psi(x,t) = (a+|\wt{\n}^{m-1} g|^2) |\wt{\n}^m g|^2, \]
where $a > 0$ is a constant to be chosen later. Exactly as before in we proceed in
the following three steps:

\begin{enumerate}
	\item Step 1: Derive an evolution inequality for $\psi$. 
	\smallskip
	
	\item Step 2: Estimate $\wt{\n} \wt{\n} \xi$ from below.
	\smallskip
	
	\item Step 3: Estimate $\xi \psi$ from above and conclude the proof.
	\medskip
	
\end{enumerate}

\noindent \textbf{Step 1: Derive an evolution inequality for $\psi$.}
\medskip

From \eqref{est nab^m g proof} and \eqref{est nab^m-1 g proof} we obtain
\begin{align}
\begin{split}
& \left(\frac{\p}{\p t} - g^{\a\b} \wt{\n}_\a \wt{\n}_\b \right)\psi \\
& = \left(\frac{\p}{\p t} - g^{\a\b} \wt{\n}_\a \wt{\n}_\b \right) ( a+|\wt{\n}^{m-1} g|^2) \cdot |\wt{\n}^m g|^2 \\
& \hspace{0.5cm}+ (a + |\wt{\n}^{m-1} g|^2) \left(\frac{\p}{\p t} - g^{\a\b} \wt{\n}_\a \wt{\n}_\b \right) |\wt{\n}^m g|^2 \\
& \hspace{0.5cm}- 2 g^{\a\b} \wt{\n}_\a |\wt{\n}^{m-1} g|^2 \wt{\n}_\b |\wt{\n}^m g|^2 \\
& \le -\frac{1}{2} |\wt{\n}^m g|^4 + c(n,m-1) |\wt{\n}^m g|^2 \cdot \\
& \hspace{0.5cm} \cdot [ C_{m-1}^2 (C_1^2 + C_2 + \sqrt{k_0}) + C_{m-1} (\mathcal{K}_{m-1} + \mathcal{L}_{m-1} + c_{m-1}) ] \\
& \hspace{0.5cm} + (a+|\wt{\n}^{m-1} g|^2) (-\frac{1}{2} |\wt{\n}^{m+1} g|^2 + c(n,m) [|\wt{\n}^m g|^2 (C_1^2 + C_2 + \sqrt{k_0}) \\ 
& \hspace{0.5cm} + |\wt{\n}^m g| (\mathcal{K}_m + \mathcal{L}_m + c_m) ]) \\
& \hspace{0.5cm} - 2 g^{\a\b} \wt{\n}_\a |\wt{\n}^{m-1} g|^2 \wt{\n}_\b |\wt{\n}^m g|^2.
\end{split}
\end{align}

We estimate the last term on the right-hand side
\begin{align}
\begin{split}
& - 2 g^{\a\b} \wt{\n}_\a |\wt{\n}^{m-1} g|^2 \wt{\n}_\b |\wt{\n}^m g|^2 \\
& \le 16 |\wt{\n}^{m-1} g| |\wt{\n}^m g|^2 |\wt{\n}^{m+1} g| \\
& \le 16 C_{m-1} |\wt{\n}^m g|^2 |\wt{\n}^{m+1} g| \\
& \le \frac{1}{2} a |\wt{\n}^{m+1} g|^2 + \frac{1}{2a} \cdot 256 C_{m-1}^2 |\wt{\n}^m g|^4.
\end{split}
\end{align}
Now choosing $a := 512 C_{m-1}^2$ yields
\begin{align}
\begin{split}
& \left(\frac{\p}{\p t} - g^{\a\b} \wt{\n}_\a \wt{\n}_\b \right)\psi \\
& \le -\frac{1}{4} |\wt{\n}^m g|^4 + c(n,m-1) |\wt{\n}^m g|^2 \cdot \\
& \hspace{0.5cm}\cdot [ C_{m-1}^2 (C_1^2 + C_2 + \sqrt{k_0}) + C_{m-1} (\mathcal{K}_{m-1} + \mathcal{L}_{m-1} + c_{m-1}) ] \\
& \hspace{0.5cm} + c(n,m) C_{m-1}^2 [|\wt{\n}^m g|^2 (C_1^2 + C_2 + \sqrt{k_0}) + |\wt{\n}^m g| (\mathcal{K}_{m} + \mathcal{L}_{m} + c_{m})].
\end{split}
\end{align}
Since 
\[ |\wt{\n}^m g|^2 = \frac{\psi}{a+ |\wt{\n}^{m-1} g|^2} \le \frac{\psi}{a} = \frac{\psi}{512 C_{m-1}^2} \]
and 
\[ |\wt{\n}^m g|^2 = \frac{\psi}{a+ |\wt{\n}^{m-1} g|^2}  \ge \frac{\psi}{(512+1) C_{m-1}^2}\]
it follows that
\begin{align}\label{psi est}
\begin{split}
& \left(\frac{\p}{\p t} - g^{\a\b} \wt{\n}_\a \wt{\n}_\b \right)\psi \\
& \le -\frac{1}{4} \frac{\psi^2}{513^2 C_{m-1}^4} + c(n,m-1) [C_1^2 + C_2 + \sqrt{k_0} + \frac{\mathcal{K}_{m-1} + \mathcal{L}_{m-1} + c_{m-1}}{C_{m-1}} ] \psi \\
& \hspace{0.5cm} + c(n,m) (C_1^2 + C_2 + \sqrt{k_0}) \psi \\
& \hspace{0.5cm} + c(n,m) C_{m-1} (\mathcal{K}_{m} + \mathcal{L}_{m} + c_{m}) \sqrt{\psi}.
\end{split}
\end{align}

\noindent \textbf{Step 2: Estimate $\wt{\n} \wt{\n} \xi$ from below.}
\medskip

\noindent Let $\eta \in C^\infty(\r)$ be the cutoff function as before and define $\xi \in C(M)$ to be the cutoff function   
\begin{align}
\begin{split}
\xi(x) = \eta \left(\frac{d_{\wt{g}}(x,x_0) - \left(\g+\frac{\d}{m+1} \right)}{\d \cdot \left(\frac{1}{2} \left( \frac{1}{m+1} + \frac{1}{m} \right) - \frac{1}{m+1} \right)}\right),
\end{split}
\end{align}
where $d_{\wt{g}}$ denotes the distance function with respect to the metric $\wt{g}$. Then 
\begin{align}
\begin{split}
& \xi(x) = 1, \qquad x \in B(x_0,\g+\frac{\d}{m+1}), \\
& \xi(x) = 0, \qquad x \in B(x_0,\g+\d) \backslash B(x_0,\g+\d \cdot \frac{1}{2} \left( \frac{1}{m+1} + \frac{1}{m} \right)), \\
& 0 \le \xi(x) \le 1, \qquad x \in M.
\end{split}
\end{align} 
If $d_{\wt{g}}(\cdot,x_0)$ is smooth in a neighborhood of a point $x$, we obtain by estimates analogous to \eqref{inbetween 3} 
\begin{align}\label{xi est 1}
|\wt{\n} \xi(x)|^2 \le \frac{c(m)}{\d^2} \xi(x)
\end{align}
and 
\begin{align}\label{xi est 2}
\begin{split}
\wt{\n}_\a \wt{\n}_\b \xi(x) \ge -\left( \frac{c(m)}{\d^2} + \frac{c(m)}{\d} \sqrt[4]{k_0} \coth \left(\sqrt[4]{k_0} d_{\wt{g}}(x,x_0) \right) \right) \wt{g}_{\a\b}(x).
\end{split}
\end{align} \ \medskip

\noindent \textbf{Step 3: Estimate $\xi \psi$ from above and conclude the proof.}
\medskip

\noindent Let
\[ F(x,t) = \xi(x)\psi(x,t), \qquad (x,t) \in B(x_0,\g+\d) \times [0,T]. \]
Since $|\wt{\n}^m g|^2(x,0) = 0$, we have
\begin{align}
F(x,0) = 0, \qquad x \in B(x_0,\g+\d).
\end{align}
Since $\xi(x) = 0$ for $x \in B(x_0,\g+\d) \backslash B(x_0,\g+\d \cdot \frac{1}{2} \left( \frac{1}{m+1} + \frac{1}{m} \right))$, it follows that
\begin{align*}
F(x,t) = 0, \qquad (x,t) \in B(x_0,\g+\d) \backslash B \left(x_0,\g+\d \cdot \frac{1}{2} \left( \frac{1}{m+1} + \frac{1}{m} \right) \right) \times [0,T].
\end{align*}
Thus there exists a point $(y_0,t_0) \in B(x_0,\g+\d \cdot \frac{1}{2} \left( \frac{1}{m+1} + \frac{1}{m} \right)) \times [0,T]$ with $t_0 > 0$ such that
\begin{align}
F(y_0,t_0) =  \max \ \left\{F(x,t) \mid (x,t) \in B(x_0,\g+\d) \times [0,T] \right\}
\end{align}
unless $F \equiv 0$ on $B(x_0,\g+\d) \times [0,T]$. \medskip

Now as before in Lemma \ref{est nab g}, 
we distinguish three cases, first case where $\xi \equiv 1$ in a neighborhood of $y_0$, 
second case where $\xi$ is not identically $1$, but smooth in a neighborhood of $y_0$, 
and third case, where $\xi$ is not smooth and a trick needs to be applied. \medskip

\noindent \underline{\textbf{Case 1. $y_0 \in B(x_0,\g+\d \cdot \frac{1}{2} \left( \frac{1}{m+1} + \frac{1}{m} \right))$}} \medskip

\noindent Then $\xi \equiv 1$ in a neighborhood of $y_0$, such that $F = \psi$ near $(y_0,t_0)$, so that
\begin{align}
\begin{split}
& 0 \le \left(\frac{\p}{\p t} - g^{\a\b} \wt{\n}_\a \wt{\n}_\b \right)\psi \\
& \le -\frac{1}{4} \frac{\psi^2}{513^2 C_{m-1}^4} + c(n,m-1) [C_1^2 + C_2 + \sqrt{k_0} + \frac{\mathcal{K}_{m-1} + \mathcal{L}_{m-1} + c_{m-1}}{C_{m-1}} ] \psi \\
& \hspace{0.5cm} + c(n,m) (C_1^2 + C_2 + \sqrt{k_0}) \psi \\
& \hspace{0.5cm} + c(n,m) C_{m-1} (\mathcal{K}_{m} + \mathcal{L}_{m} + c_{m}) \sqrt{\psi}
\end{split}
\end{align}
and hence 
\begin{align}
\begin{split}
\frac{1}{4} \frac{F^2(y_0,t_0)}{513^2 C_{m-1}^4} \le & \; c(n,m-1) [C_1^2 + C_2 + \sqrt{k_0} + \frac{\mathcal{K}_{m-1} + \mathcal{L}_{m-1} + c_{m-1}}{C_{m-1}} ] F(y_0,t_0) \\
& + c(n,m) (C_1^2 + C_2 + \sqrt{k_0}) F(y_0,t_0) \\
& + c(n,m) C_{m-1} (\mathcal{K}_{m} + \mathcal{L}_{m} + c_{m}) \sqrt{F(y_0,t_0)},
\end{split}
\end{align}
which again is a better estimate than the one below in Case 2, and hence Case 1 follows from Case 2. \\

\noindent \underline{\textbf{Case 2. $y_0 \notin B(x_0,\g+\d \cdot \frac{1}{2} \left( \frac{1}{m+1} + \frac{1}{m} \right))$ and $y_0$ is not in the cut locus of $x_0$}}\medskip

\noindent Then the distance function $d_{\wt{g}}(\cdot,x_0)$ is smooth in a neighborhood of $y_0$ and we have
\begin{align}\label{F inequ at max}
\begin{split}
& 0 \le \frac{\p F}{\p t}(y_0,t_0) = \xi(y_0) \frac{\p \psi}{\p t}(y_0,t_0), \\
& 0 = \wt{\n}_\a F(y_0,t_0) = (\xi \wt{\n}_\a \psi + \psi \wt{\n}_\a \xi)(y_0,t_0), \\
& 0 \ge g^{\a\b} \wt{\n}_\a \wt{\n}_\b F(y_0,t_0) = (\xi g^{\a\b} \wt{\n}_\a \wt{\n}_\b \psi + \psi g^{\a\b} \wt{\n}_\a \wt{\n}_\b \xi 
\\ &\hspace{0.7cm}+ 2 g^{\a\b} \wt{\n}_\a \xi \wt{\n}_\b \psi)(y_0,t_0).
\end{split}
\end{align}
Together with \eqref{psi est} we obtain at the point $(y_0,t_0)$
\begin{align} \label{psi xi est}
\begin{split}
\frac{1}{4}\frac{\psi^2}{513^2 C_{m-1}^4}\xi & \le \xi g^{\a\b}\wt{\n}_\a \wt{\n}_\b \psi + \xi \cdot \mathcal{A}\\
& \le -\psi g^{\a\b} \wt{\n}_\a \wt{\n}_\b \xi - 2 g^{\a\b} \wt{\n}_\a \xi \wt{\n}_\b \psi + \xi \cdot \mathcal{A}
\end{split}
\end{align}
with
\begin{align}
\begin{split}
 \mathcal{A} & := c(n,m-1) [C_1^2 + C_2 + \sqrt{k_0} + \frac{\mathcal{K}_{m-1} + \mathcal{L}_{m-1} + c_{m-1}}{C_{m-1}} ] \psi \\
& \hspace{0.5cm}+ c(n,m) (C_1^2 + C_2 + \sqrt{k_0}) \psi \\
& \hspace{0.5cm}+ c(n,m) C_{m-1} (\mathcal{K}_{m} + \mathcal{L}_{m} + c_{m}) \sqrt{\psi}.
\end{split}
\end{align}
Using \eqref{F inequ at max}, \eqref{xi est 1} and \eqref{equiv met 2} we have at $(y_0,t_0)$
\begin{align} \label{est 1}
- 2 g^{\a\b} \wt{\n}_\a \xi \wt{\n}_\b \psi  = \frac{2\psi}{\xi} g^{\a\b} \wt{\n}_\a \xi \wt{\n}_\b \xi \le \frac{c(m)}{\d^2} \psi.
\end{align}
Also \eqref{xi est 2} and an estimate analogous to \eqref{est coth} yields at $(y_0,t_0)$
\begin{align} \label{est 2}
-\psi g^{\a\b} \wt{\n}_\a \wt{\n}_\b \xi \le \frac{c(n,m,k_0)}{\d^2} \psi.
\end{align}
Plugging \eqref{est 1} and \eqref{est 2} into \eqref{psi xi est} leads to
\begin{align}
\begin{split}
\frac{1}{4}\frac{\psi^2}{513^2 C_{m-1}^4}\xi \le & \frac{c(n,m,k_0)}{\d^2} \psi + \xi \cdot \mathcal{A}. 
\end{split}
\end{align}
Multiplying by $\xi$ while using $0 \le \xi \le 1$ and adjusting the constants $c(n,m,k_0)>0$ we obtain
\begin{align}\label{second-derivatie-estimate1}
\begin{split}
\frac{F(y_0,t_0)^2}{C_{m-1}^4} & \le c(n,m-1) [C_1^2 + C_2 + \sqrt{k_0} + \frac{\mathcal{K}_{m-1} + \mathcal{L}_{m-1} + c_{m-1}}{C_{m-1}} ] F(y_0,t_0) \\
& \hspace{0.5cm}+ c(n,m) (C_1^2 + C_2 + \sqrt{k_0}) F(y_0,t_0) \\
& \hspace{0.5cm}+ c(n,m) C_{m-1} (\mathcal{K}_{m} + \mathcal{L}_{m} + c_{m}) \sqrt{F(y_0,t_0)}
\\ &\hspace{0.5cm}+ \frac{c(n,m,k_0)}{\d^2} F(y_0,t_0).
\end{split}
\end{align}
Now we use the elementary estimate \eqref{elem est} with $a = 0$, namely: If $x \ge 0$ satisfies
\[ x^2 \le bx + cx^{1/2}\]
with constants $b,c \ge 0$, then 
\[ x \le \max\{b,c^{2/3}\}.\] 
This reduces \eqref{second-derivatie-estimate1} to 
\begin{align}
\begin{split}
F(y_0,t_0) & \le \max \{ c(n,m,k_0) C_{m-1}^4 \left( \frac{1}{\d^2} + C_1^2 + C_2 + \sqrt{k_0} + \frac{\mathcal{K}_{m-1} + \mathcal{L}_{m-1} + c_{m-1}}{C_{m-1}}    \right), \\
& \hspace{1.7cm}c(n,m) (C_{m-1}^5 (\mathcal{K}_{m} + \mathcal{L}_{m} + c_{m}))^{2/3}   \}.
\end{split}
\end{align}
Hence for all $(x,t) \in B(x_0,\g+\d) \times [0,T]$
\begin{align}
\begin{split}
& F(x,t) \le F(y_0,t_0) \le \\
& \max \{ c(n,m,k_0) C_{m-1}^4 \left( \frac{1}{\d^2} + C_1^2 + C_2 + \sqrt{k_0} + \frac{\mathcal{K}_{m-1} + \mathcal{L}_{m-1} + c_{m-1}}{C_{m-1}}    \right), \\
& \hspace{1.2cm} c(n,m) (C_{m-1}^5 (\mathcal{K}_{m} + \mathcal{L}_{m} + c_{m}))^{2/3}   \}.
\end{split}
\end{align}
As $F(x,t) = \xi(x) \psi(x,t)$ and $\xi(x) = 1$ for $x \in B(x_0,\g+\frac{\d}{m+1})$, we conclude
\begin{align}
\begin{split}
\psi(x,t) & \le \max \{ c(n,m,k_0) C_{m-1}^4 \left( \frac{1}{\d^2} + C_1^2 + C_2 + \sqrt{k_0} + \frac{\mathcal{K}_{m-1} + \mathcal{L}_{m-1} + c_{m-1}}{C_{m-1}}    \right), \\
& \hspace{1.7cm} c(n,m) (C_{m-1}^5 (\mathcal{K}_{m} + \mathcal{L}_{m} + c_{m}))^{2/3}   \}
\end{split}
\end{align}
for all $(x,t) \in B(x_0,\g+\frac{\d}{m+1}) \times [0,T]$.
Now since $\psi(x,t) = (a + |\wt{\n}^{m-1} g|^2) |\wt{\n}^m g|^2$ and $a = 512 C_{m-1}^2$ we finally obtain
\begin{align} 
\begin{split}
& |\wt{\n}^m g|^2 (x,t) = \frac{\psi(x,t)}{(a + |\wt{\n}g|^2)(x,t)} \le \frac{\psi(x,t)}{a} \\
& \le \max \{ c(n,m,k_0) C_{m-1}^2 \left( \frac{1}{\d^2} + C_1^2 + C_2 + \sqrt{k_0} + \frac{\mathcal{K}_{m-1} + \mathcal{L}_{m-1} + c_{m-1}}{C_{m-1}}    \right), \\
& \hspace{1.7cm}c(n,m) \frac{1}{C_{m-1}^2}(C_{m-1}^5 (\mathcal{K}_{m} + \mathcal{L}_{m} + c_{m}))^{2/3}   \} 
\end{split}
\end{align}
for all $(x,t) \in B(x_0,\g+\frac{\d}{m+1}) \times [0,T]$. \medskip

\noindent \underline{\textbf{Case 3. $y_0 \notin B(x_0,\g+\d \cdot \frac{1}{2} \left( \frac{1}{m+1} + \frac{1}{m} \right))$ and $y_0$ is in the cut locus of $x_0$}} \medskip

\noindent Here we again apply Calabi's trick, see Case 3 in the Proof of \cref{est nab g}.

\end{proof}

\section{Proof of the main existence and regularity result}\label{main-section}

In this section we describe the necessary modifications of the proofs in Section 
\ref{Review of Shis local existence theorem} for the case of incomplete manifolds.
Despite incompleteness, we still continue under the assumption of bounded geometry
$|\wt{\Rm}|^2 \le k$ for some positive constant $k>0$. 

\subsection{Validity of Theorems \ref{Shi apriori} and \ref{Shi Theorem 3.2} in the incomplete case} \ \medskip

We start by observing that, due to the relative compactness of the domain $D$, 
\cref{Shi apriori} and \cref{Shi Theorem 3.2} are also valid in case the initial manifold $(M,\wt{g})$ is 
incomplete. Indeed, we can follow the same steps as in the proof of \cref{Shi apriori} outlined above. 
To give just one example from the proof of \cref{Shi Theorem 3.2}, the compactness of the closure of the 
domain $D$ still gives a positive lower bound for the injectivity radius on $\ov{D}$, which is needed for 
the estimate (7) in Lemma 3.1 in \cite{shi}. \medskip

Notice that  the injectivity radius on $\ov{D}$ can become small in the case of incomplete manifolds when $D$ is close to the singularity, but that smallness of the injectivity radius can also happen in the  complete case, e.g. on manifolds with  hyperbolic cusps. At this point we emphasize once again. cf. Remark \ref{lower-injectivity}, that the lower bound for the injectivity radius on $\ov{D}$ does not enter the existence times $T(n,k_0)$ and $T(n,k_0, \delta)$. This is a crucial point, as we will later take an exhaustion of the manifold $M$ by such domains which then get closer and closer to the singularity.

\subsection{Extension of Lemmas \ref{Shi Lemma 4.1} and \ref{Shi Lemma 4.2} to the incomplete case} \ \medskip

Next we formulate and prove interior estimates for the derivatives of the metric, corresponding to 
\cref{Shi Lemma 4.1} and \cref{Shi Lemma 4.2}. Let $U \subset M$ be open and relatively compact, such that $\p U$ is an $(n-1)$-dimensional, smooth, compact submanifold. Choose $\delta > 0$ small enough, that $\overline{B(\ov{U},\delta)} \subset M$ is compact and that the function $d_{\wt{g}}(\cdot,\ov{U}): M \to \r$ giving the distance to $\ov{U}$ is smooth on $B(\ov{U},\delta) \backslash \ov{U}$. The latter is possible, since $d_{\wt{g}}(x,\ov{U}) = d_{\wt{g}}(x,\p U)$ for all $x \in M \backslash \ov{U}$, and $d_{\wt{g}}(\cdot,\p U)$ is smooth in a neighborhood of $\p U$ by \cite[Theorem 1  and Remark (1)]{fo}. The following result is an extension of Lemmas \ref{Shi Lemma 4.1} and \ref{Shi Lemma 4.2} to the incomplete case.

\begin{lem}\label{Lemma 4.2}
Fix $U, \d$ as above, and a finite $T>0$. Let $g(x,t)$ be a smooth solution of the initial value problem
\begin{equation*}
\begin{split}
\frac{\p}{\p t} g_{ij}(x,t) = (-2 \Ric_{ij} + \n_i V_j + \n_j V_i)(x,t)&, \quad  (x,t) \in \, B(\ov{U},\d) \times [0,T],  \\
g(x,0) = \wt{g}(x)&, \quad x \in \, B(\ov{U},\d),
\end{split}
\end{equation*}
where $V^i = g^{jk}(\G^i_{jk} - \wt{\G}^i_{jk})$ is the de Turck vector field. Furthermore, assume that
\begin{equation*}
(1-\ve(n)) \wt{g}(x) \le g(x,t) \le (1+\ve(n)) \wt{g}(x)
\end{equation*}
for $\ve(n) > 0$ sufficiently small only depending on $n$ and for all $(x,t) \in B(\ov{U},\d) \times [0,T]$. Then for all $m \in \N_0$ there exists $c(n,m,U,\d,T,\wt{g}) > 0$ depending only on $n, m, U,\d,T$ and $\wt{g}$, such that for all $(x,t) \in B(\ov{U},\frac{\d}{m+1}) \times [0,T]$.
\begin{equation*}
|\wt{\n}^m g(x,t)|^2 \le c(n,m,U,\d,T,\wt{g}).
\end{equation*}
\end{lem}

\begin{proof}
Consider exactly as in \eqref{phi} the function
\[ \varphi(x,t) := a + \sum_{k=1}^n \l_k(x,t)^{m_0}, \quad (x,t) \in B(\ov{U}, \d) \times [0,T], \]
where $a,m_0$ are the same constants as in the proof of Lemma \ref{Shi Lemma 4.1} only depending on $n$,  
and $\l_k(x,t)$ are the eigenvalues of $g(x,t)$ with respect to $\wt{g}(x)$. Let
\[ \psi(x,t) := |\wt{\n}g|^2 \varphi(x,t). \]
These are the same functions as in \eqref{phi}, \eqref{psi}, but with $B(x_0,\g + \eta)$ 
replaced by $B(\ov{U},\d)$. Then performing the same steps as in the proof of \cref{Shi Lemma 4.1}, we obtain
\[ \frac{\p \psi}{\p t} \le g^{\a\b} \wt{\n}_\a \wt{\n}_\b \psi - \frac{1}{16} \psi^2 + c_0, \]
on $B(\ov{U}, \d) \times [0,T]$, where $c_0 > 0$ is a constant only depending on $n$ and $\wt{g}$. \medskip

Now let $\eta \in C^\infty(\r)$ be a smooth, nonincreasing cutoff function, 
such that $\eta \equiv 1$ on $(-\infty, 0]$, $\eta \equiv 0$ on $[1,\infty)$, satisfying \eqref{eta properties}. 
Such a function is illustrated in Figure \ref{fig:CutOff1} above. Then we define a cutoff function $\xi \in C^\infty_c(M)$ via
\[ 
\xi(x) := \eta \left(\frac{d_{\wt{g}}(x,\ov{U}) - \d/2}{\d/4}\right), \quad 
\textup{for all} \ x \in M, 
\]
which is a modified version of \eqref{xi}. Observe that $\xi$ satisfies
\begin{equation}
\begin{split}
  \xi(x) = 1, & \quad x \in B(\ov{U},\d/2), \\
  \xi(x) = 0, & \quad x \in M \backslash B(\ov{U},3 \d/4).
 \end{split}
\end{equation}
Exactly as before in \eqref{prop cutoff 1}, we can control the gradient of $\xi$
\begin{equation}
|\wt{\n} \xi|^2 (x) \le \frac{16^2}{\d^2} \xi(x), \quad x \in M,
\end{equation}
since the distance function $d_{\wt{g}}(\cdot,\ov{U})$ is smooth with $|\wt{\n} d_{\wt{g}}(x,\ov{U})| \le 1$
on $B(\ov{U},\delta) \backslash \ov{U}$. Since $\ov{B(\ov{U},3 \d/4)} \subset M$ is compact, 
we still find exactly as in \eqref{prop cutoff 2}
\begin{equation}
\wt{\n} \wt{\n} \xi \ge -c\wt{g}
\end{equation}
for some constant $c > 0$ only depending on $U, \d$ and $\wt{g}$.
Define
\[ F(x,t) := \xi(x) \psi(x,t), \quad (x,t) \in B(\ov{U}, \d) \times [0,T]. \]
By construction, it has the properties
\begin{equation}
\begin{split}
&F(x,0) = 0, \quad x \in B(\ov{U}, \d), \\
&F(x,t) = 0, \quad (x,t) \in M \backslash B(\ov{U}, 3\d/4) \times [0,T].
\end{split}
\end{equation}
Hence there exists $(x_0,t_0) \in B(\ov{U}, 3\d/4) \times [0,T]$ with
\[ 
F(x_0,t_0) = \max \, \{ \, F(x,t) \mid (x,t) \in B(\ov{U}, 3\d/4) \times [0,T] \}. 
\]
Now following the same steps as in the proof of \cref{Shi Lemma 4.1} we obtain
\[ F(x_0,t_0) \le c(n,U,\d,T,\wt{g}), \]
for some constant $c(n,U,\d,T,\wt{g})>0$ depending only on the arguments $n,U,\d,T$
and $\wt{g}$. We conclude
\[ 
\xi(x) \psi(x,t) \equiv F(x,t) \le c(n,U,\d,T,\wt{g}), \quad (x,t) \in B(\ov{U}, \d) \times [0,T]. 
\]
Since $ \xi \equiv 1$ on $B(\ov{U},\d/2)$, we find
\[ 
|\wt{\n}g|^2 \varphi(x,t) \equiv \psi(x,t) \le c(n,U,\d,T,\wt{g}), \quad (x,t) \in B(\ov{U}, \d/2) \times [0,T]. 
\]
The claim now follows in case $m=1$ once we divide both sides of the inequality by $\varphi(x,t) \ge a > 0$,
which gives
\[ 
|\wt{\n}g|^2(x,t) \le \frac{1}{a}c(n,U,\d,T,\wt{g}), \quad (x,t) \in B(\ov{U}, \d/2) \times [0,T]. 
\]
The cases $m \ge 2$ are proven by induction. Set as in the proof of \cref{Shi Lemma 4.2}  
\[ \Psi(x,t) := (a_0+ |\wt{\n}^{m-1} g(x,t)|^2) |\wt{\n}^m g(x,t)|^2. \]
Replacing balls $B(x_0,\g + \d/k)$ by neighborhoods $B(\ov{U},\d/k)$, and performing the same steps as in \cite[Proof of Lemma 4.2]{shi}, meanwhile choosing $a_0 > 0$ only depending on $m,n,U,\d,T, \wt{g}$ appropriately, we obtain
\[ \frac{\p \Psi}{\p t} \le g^{\a\b} \wt{\n}_\a \wt{\n}_\b \Psi - c_1 \Psi^2 + c_0, \]
on $B(\ov{U},\d/m) \times [0,T]$, with constants $c_0, c_1 > 0$ only depending on $m,n,U,\d,T$ and $\wt{g}$. 
Now using the cutoff function
\[ 
\xi_m (x) := \eta \left(\frac{d_{\wt{g}}(x,\ov{U}) - \frac{\d}{m+1}}{\frac{\d}{m+\frac{1}{2}} - \frac{\d}{m+1}}\right), \quad x \in M, 
\]
we obtain similar to the proof of \cref{Shi Lemma 4.2}
\[ 
\Psi(x,t) \le c_2(m,n,U,\d,T,\wt{g}), \quad (x,t) \in B(\ov{U},\d/(m+1)) \times [0,T]. 
\]
This leads to the inequality
\[ 
|\wt{\n}^m g(x,t)|^2 \le \frac{1}{a_0} \Psi(x,t) \le \frac{1}{a_0} c_2(m,n,U,\d,T,\wt{g}), 
\] 
for all $(x,t) \in B(\ov{U},\d/(m+1)) \times [0,T]$,
which proves the case $m \ge 2$.
\end{proof}

\subsection{Proof of \cref{main theo} by exhaustion} \ \medskip

\noindent Now we can finish the proof of \cref{main theo}. Let $\{\ov{U_k}\}_{k \in \N}$, be an exhaustion of $M$ by $n$-dimensional, smooth, compact manifolds with boundary, i.e. $U_k \subset M$ is open, $\ov{U_k} \subset M$ is compact, $\p U_k$ is an $(n-1)$-dimensional, smooth, compact submanifold, $U_k \subset U_{k+1}$ for all $k \in \N$ and $\bigcup\limits_{k \in \N} U_k = M$. By \cref{Shi apriori} and \cref{Shi Theorem 3.2} there exists $T(n,k_0) > 0$, such that the system (cf. \eqref{Dirichlet Ricci-k})
\begin{equation}\label{BVP}
\begin{split}
\frac{\p}{\p t} g_{ij}(x,t) = (-2 \Ric_{ij} + \n_i V_j + \n_j V_i)(x,t)&, \quad  (x,t) \in U_k \times [0,T(n,k_0)],  \\
g(x,t) = \wt{g}(x)&, \quad (x,t) \in \partial U_k \times [0,T(n,k_0)], \\
g(x,0) = \wt{g}(x)&, \quad x \in U_k.
\end{split}
\end{equation}
where $V^i = g^{jk}(\G^i_{jk} - \wt{\G}^i_{jk})$ is the de Turck vector field, has a unique smooth solution $g(k,x,t)$ on $0 \le t \le T(n,k_0)$ for all $k \in \N$, satisfying the estimate
\begin{equation}
(1-\ve(n)) \wt{g}(x) \le g(k,x,t) \le (1+\ve(n)) \wt{g}(x)
\end{equation}
for all $(x,t) \in \ov{U_k} \times [0,T(n,k_0)]$ and for $\ve(n) > 0$ sufficiently small only depending on $n$. Note that $\ov{U_k}$ from the exhaustion above need not be connected, but since it is compact, it has at most finitely many connected components, so that \cref{Shi apriori} and \cref{Shi Theorem 3.2} can be applied to each component. \medskip

Choose $\d_k > 0$ sufficiently small, such that the closure of $B(\ov{U_k},\d_k) \subset M$ is compact and 
such that the function $d_{\wt{g}}(\cdot,\ov{U_k}): M \to \r$ is smooth on $B(\ov{U}_k,\d_k) \backslash \ov{U_k}$.
By compactness of the closure of $B(\ov{U_1},\d_1)$, there exists $N \in \N$, such that the solution $g(k,x,t)$ is 
defined on $B(\ov{U_1},\d_1)$ for all $k \ge N$. By \cref{Lemma 4.2}
\[ 
|\wt{\n}^m g(k,x,t)|^2 \le c(n,m,U_k,\d_k,T,\wt{g}), 
\]
for all $k \ge N$, $m \in \N_0$, and $(x,t) \in \ov{U_1} \times [0,T(n,k_0)]$. 
Then by Arzel\`a-Ascoli there exists a subsequence $(g(k_\ell,x,t))_{\ell \in \N}$, which 
converges on $\ov{U_1} \times [0,T]$ in the $C^\infty$ topology to a 
family of $C^\infty$ metrics $g(x,t)$. \medskip

Similarly a subsequence of the subsequence converges on $\ov{U_2} \times [0,T]$, etc. Now the diagonal sequence converges on every $\ov{U_k} \times [0,T]$ to $g(x,t)$. As the sequence $(U_k)$ eventually contains any given compact subset of $M$, the diagonal sequence converges smoothly locally uniformly to $g(x,t)$. Then $g(x,t)$ solves \eqref{Ricci de Turck system on M incomplete}. The estimate \eqref{unif equiv delta on M incomplete} follows by restricting the solutions $g(k,x,t)$ to $t \in [0,T(n,k_0,\d)]$, where $T(n,k_0,\d)$ is from \cref{Shi apriori}. 

\section{Open problems and future research directions} 

We intend to discuss the following questions in the subsequent publications.

\begin{enumerate}
\item Does the Ricci de Turck flow, presented here, and the flow constructed 
by the second author in \cite{v}, coincide in the setting of incomplete manifolds of bounded
geometry with wedge singularities?
\item Can we extend the tensor maximum principle to the incomplete setting?
\item Is there a way to define a flow of arbitrary incomplete manifolds without assuming bounded curvature, 
for instance imposing bounded Ricci curvature only? 
\end{enumerate}


\begin{thebibliography}{AAAA}

\setlength{\parsep}{0ex}

\bibitem[\textsc{BaVe14}]{BV}
Eric Bahuaud and Boris Vertman, \emph{Yamabe flow on manifolds with edges},
Math. Nachr. \textbf{287} (2014), no.~23, 127--159 

\bibitem[\textsc{CLN06}]{CLN}
B.~Chow, P.~Lu, L.~Ni, \textit{Hamilton's Ricci flow}, AMS Science Press, 2006

\bibitem[\textsc{Foo84}]{fo}
R.~Foote, \textit{Regularity of the distance function},
Proc. of the AMS \textbf{92}, no. 1 (1984), 153--155

\bibitem[\textsc{GiTo11}]{Topping2}
\bysame, \emph{Existence of {R}icci flows of incomplete surfaces}, Comm.
Partial Differential Equations \textbf{36} (2011), no.~10, 1860--1880

\bibitem[\textsc{Ham82}]{ham}
R.~Hamilton, \textit{Three-manifolds with positive Ricci curvature},
J. Differential Geom. \textbf{17} (1982), 255--306

	
\bibitem[\textsc{KrVe19a}]{Klaus-Vertman}
Klaus Kr\"{o}ncke and Boris Vertman,
\newblock \emph{Stability of Ricci de Turck flow on Singular Spaces},
{Calc. Var. Partial Differ. Equ.} \textbf{58} (2019),  no. 2, 74 

\bibitem[\textsc{MRS15}]{MRS}
Rafe Mazzeo, Yanir Rubinstein, and Natasha Sesum, \emph{Ricci flow on surfaces with conic singularities}, 
Anal. PDE \textbf{8} (2015), no.~4, 839-–882

\bibitem[\textsc{Per02}]{per1}
G.~Perelman, \textit{The entropy formula for the Ricci flow and its geometric applications},
arXiv:math/0211159v1 [math.DG]

\bibitem[\textsc{Per03a}]{per2}
G.~Perelman, \textit{Ricci flow with surgery on three-manifolds},
arXiv:math/0303109v1 [math.DG]

\bibitem[\textsc{Per03b}]{per3}
G.~Perelman, \textit{Finite extinction time for the solutions to the Ricci flow on certain three-manifolds},
arXiv:math/0307245v1 [math.DG]

\bibitem[\textsc{Pet06}]{pet}
Peter Petersen, \textit{Riemannian Geometry}, Second Edition, Springer, 2006

\bibitem[\textsc{Shi89}]{shi}
W.-X.~Shi, \textit{Deforming the metric on complete Riemannian manifolds},
J. Differential Geometry \textbf{30} (1989), 223--301



\bibitem[\textsc{Sim13}]{MS}
Miles Simon, \emph{Local smoothing results for the Ricci flow in dimensions two and three}, 
Geom. Topol. \textbf{17} (2013), no. 4, 2263--2287


\bibitem[\textsc{Ver16}]{v}
B.~Vertman, \textit{Ricci flow on singular manifolds},\\
arXiv:1603.06545v3 [math.DG] (2016)

\bibitem[\textsc{Yin10}]{Yin:RFO}
Hao Yin, \emph{Ricci flow on surfaces with conical singularities}, J. Geom.
Anal. \textbf{20} (2010), no.~4, 970--995.


\end{thebibliography}
\end{document}